\documentclass[12pt]{amsart}
\usepackage{amssymb,amsmath,amsthm}

\numberwithin{equation}{section}
\begin{document}

\theoremstyle{plain}
\newtheorem{theorem}{Theorem}[section]
\newtheorem{lemma}[theorem]{Lemma}
\newtheorem{proposition}[theorem]{Proposition}
\newtheorem{corollary}[theorem]{Corollary}
\newtheorem{conjecture}[theorem]{Conjecture}

\def\mod#1{{\ifmmode\text{\rm\ (mod~$#1$)}
\else\discretionary{}{}{\hbox{ }}\rm(mod~$#1$)\fi}}

\theoremstyle{definition}
\newtheorem*{definition}{Definition}

\theoremstyle{remark}
\newtheorem*{remark}{Remark}
\newtheorem{example}{Example}[section]
\newtheorem*{remarks}{Remarks}

\newcommand{\qq}{{\mathbb Q}}
\newcommand{\rr}{{\mathbb R}}
\newcommand{\nn}{{\mathbb N}}
\newcommand{\zz}{{\mathbb Z}}
\newcommand{\al}{\alpha}
\newcommand{\be}{\beta}
\newcommand{\ga}{\gamma}
\newcommand{\mz}{{\mathcal Z}}
\newcommand{\mi}{{\mathcal I}}
\newcommand{\ep}{\epsilon}
\newcommand{\la}{\lambda}
\newcommand{\de}{\delta}
\newcommand{\Del}{\Delta}
\title{On Hilbert's construction of positive polynomials}

\thanks{This material is based in part upon work of the author,
  supported by the USAF under DARPA/AFOSR MURI Award
  F49620-02-1-0325. Any opinions, findings, and conclusions or
  recommendations expressed in this publication are those of the
  author and do not necessarily reflect the views of these agencies.}

\author{Bruce Reznick}
\address{Department of Mathematics, University of 
Illinois at Urbana-Champaign, Urbana, IL 61801} 
\email{reznick@math.uiuc.edu}
\subjclass{Primary: 11E20, 11E25, 12D99,  14H50, 14N15}
\begin{abstract}
In 1888, Hilbert described how to find real polynomials which take
only non-negative values but are not a sum of squares of
polynomials. His construction was so restrictive that no explicit
examples appeared until the late 1960s. We revisit and
generalize Hilbert's construction and present many such polynomials.
\end{abstract}
\date{\today}

\maketitle

\section{History and Overview}
A real polynomial $f(x_1,\dots,x_n)$ is {\it psd} or {\it positive} if
$f(a) \ge 
0$ for all $a \in \mathbb R^n$; it is {\it sos} or a {\it sum of
  squares} if there exist real polynomials $h_j$ so that $f = \sum h_j^2$.
For forms, we follow the notation
of \cite{CL1} and use  $P_{n,m}$ to denote the  cone of real psd
forms of even degree $m$ in $n$ variables, $\Sigma_{n,m}$ to denote its
subcone of sos forms and let $\Delta_{n,m} = P_{n,m} \smallsetminus
\Sigma_{n,m}$.  The 
Fundamental Theorem of Algebra implies that $\Delta_{2,m} = \emptyset$;
 $\Delta_{n,2} = \emptyset$ follows from the diagonalization of psd 
quadratic forms.

The first suggestion that a psd form might not be sos was made by
Minkowski in the oral defense of his 1885 doctoral dissertation: Minkowski
proposed the thesis that not every psd form is sos. Hilbert was one of his
official ``opponents'' and remarked that Minkowski's
arguments had convinced him that this thesis should be true for
ternary forms. (See \cite{Hi3}, \cite{Min} and \cite{Sch}.) Three years later,
in a single remarkable paper, Hilbert \cite{Hi1} resolved the question.
He first showed that $F \in P_{3,4}$ is a sum of three squares of
quadratic forms; see \cite{Ru} and \cite{Sw} for recent 
expositions and \cite{PR,PRSS} for another approach.
Hilbert then described a construction of
forms in $\Delta_{3,6}$ and $\Delta_{4,4}$;
after multiplying these by powers of linear forms if necessary, it
follows that $\Delta_{n,m} \neq \emptyset$ if $n \ge 3$ and $m \ge 6$ or
$n \ge 4$ and $m \ge 4$.

 The goal of this paper is to
isolate the underlying mechanism of Hilbert's construction, show that
it applies to  situations more general than those in
\cite{Hi1}, and use it to produce many new examples.

In \cite{Hi1}, Hilbert first worked with polynomials in two variables,
which homogenize to ternary forms.
Suppose $f_1(x,y)$ and $f_2(x,y)$ are two relatively
prime real cubic polynomials with  nine distinct real
common zeros -- $\{\pi_i\}$, indexed arbitrarily -- so that no three
of the $\pi_i$'s lie on a line and no six lie on a quadratic. By counting
coefficients, one sees that there exists a non-zero quadratic $\phi(x,y)$
with zeros at 
$\{\pi_1,\dots,\pi_5\}$ and a non-zero quartic  $\psi(x,y)$ with the
same zeros, and which is singular at $\{\pi_6,\pi_7,\pi_8\}$:
 the sextic $\phi\psi$ is thus singular at $\{\pi_1,\dots,\pi_8\}$. Hilbert
showed that $(\phi\psi)(\pi_9) \neq 0$ and that there exists
 $c \neq 0$ so that the perturbed polynomial
$p = f_1^2 + f_2^2 + c\phi\psi$ is positive. If $p =
\sum h_j^2$, then each $h_j$ would be a cubic which vanishes on
$\{\pi_1,\dots,\pi_8\}$. But Cayley-Bacharach implies that $h_j(\pi_9) = 0$ for
each $j$, hence $p(\pi_9) = 0$, a contradiction. Thus, $p$ homogenizes
to a form $P \in \Delta_{3,6}$.

Hilbert also considered  in \cite{Hi1} three relatively prime real quadratic
polynomials, $f_i(x,y,z)$, $1 \le i \le 3$, with 
eight distinct real common zeros -- $\{\pi_i\}$, indexed arbitrarily --
so that no four of the zeros lie on a plane. There exists
a non-zero linear $\phi(x,y,z)$ with zeros at 
$\{\pi_1,\pi_2,\pi_3\}$ and a non-zero cubic $\psi(x,y,z)$ with the
same zeros, and which is singular at $\{\pi_4,\pi_5,\pi_6,\pi_7\}$. 
Similarly, $(\phi\psi)(\pi_8) \neq 0$ and there exists
 $c \neq 0$ so that $f_1^2 + f_2^2 + f_3^2+ c\phi\psi$ is positive and not
sos. This homogenizes to a form in $\Delta_{4,4}$.
 
In 1893, Hilbert \cite{Hi2} showed that if
$F \in P_{3,m}$ with $m \ge 4$, then there exists a form 
$G \in P_{3,m-4}$ and forms $H_{k}$, $1 \le k \le 3$, so that
$GF= H_{1}^2 + H_{2}^2 + H_{3}^2$. 
(Hilbert's construction does not readily identify $G$ or the $H_k$'s.)
In particular, if $F \in P_{3,6}$, then there exists 
$Q\in P_{3,2}$ so that $QF \in \Sigma_{3,8}$; since $Q\cdot QF \in
\Sigma_{3,10}$, $F$ is a sum of squares of rational functions with
common denominator $Q$. An iteration of this argument shows
that if $F \in P_{3,m}$, then there exists $G$ so that $G^2F$ is sos. 
Hilbert's 17th Problem \cite{Hi23} asked whether this representation
as a sum of squares of rational functions exists for forms in
$P_{n,m}$ when $n \ge 4$.
For much more on the history of this subject up to 1999, see
\cite{Re2}. Recently, Blekherman \cite{Bl}
has shown that for fixed degree $m$, the ``probability'' that a psd
form is sos goes to 0 as $n$ increases. This result highlights the importance
of understanding psd forms which are not sos.

Hilbert's restriction on the common zeros meant
that no very simple or symmetric example could be constructed, and the first
explicit example of any $P \in \Delta_{n,m}$ did not appear for many
decades. The only two detailed references to Hilbert's construction
before the late 1960s (known to the author) are by 
Terpstra \cite{Ter} (on biquadratic forms, related to
$\Delta_{4,6}$, thanks to Roland Hildebrand for the reference), and an 
exposition \cite[pp.232-235]{GV} by Gel'fand and Vilenkin of the
sextic case only. 

At a 1965 conference on inequalities,
 Motzkin \cite{Mo}  presented a specific sextic polynomial
$m(x,y)$  which is positive by the 
arithmetic-geometric inequality and not sos by the
arrangement of monomials in its Newton polytope. (Hilbert's last
assistant, Olga Taussky-Todd, who
had a lifelong interest in sums of squares, heard Motzkin speak, and
informed him 
that $m(x,y)$ was the first specific polynomial known to be positive but not
sos.)   After
homogenization, Motzkin's example is
\begin{equation}
M(x,y,z) = x^4y^2 + x^2y^4 + z^6 - 3x^2y^2z^2 \in \Delta_{3,6}.
\end{equation}
Around the same time and independently,
R. M. Robinson \cite[p.264]{Ro} wrote that he
saw ``an unpublished example of a ternary
sextic worked out recently by W. J. Ellison using Hilbert's Method. It
is, as would be expected,
very complicated. After seeing this, I discovered that an astonishing
simplification would be possible
by dropping some unnecessary assumptions made by Hilbert." 
Robinson observed that the cubics $f_1(x,y) = x^3-x$ and
$f_2(x,y) = y^3-y$ have nine common zeros: the $3 \times 3$ square
$\{-1,0,1\}^2$.  There are eight lines which each contain three of the
zeros. Still, the sextic $(x^2-1)(y^2-1)(1-x^2-y^2)$ is positive
at (0,0) and singular at the other eight points. By taking the
maximum value for $c$ in Hilbert's construction and homogenizing,
Robinson showed that 
\begin{equation}
R(x,y,z) = x^6 + y^6 + z^6 - x^4y^2 - x^2y^4 - x^4z^2-y^4z^2 -
x^2z^4-y^2z^4+3x^2y^2z^2
\end{equation}  
is in $\Delta_{3,6}$. Similarly, by taking the three quadratics $x^2-x$,
$y^2-y$ and $z^2-z$, whose common zeros are $\{0,1\}^3$,
 choosing $(1,1,1)$ as the eighth point, and then
homogenizing,  Robinson showed that
\begin{equation}
\tilde R(x,y,z,w) = x^2(x-w)^2 + y^2(y-w)^2+z^2(z-w)^2+2xyz(x+y+z-2w)
\end{equation}
is in $\Delta_{4,4}$. (The only other published implementation of
Hilbert's Method known to the author  is a 1979 sextic studied by Schm\"udgen
\cite{Schm} using $\{-2,0,2\}^2$, with ninth
point $(2,0)$.)

The papers of Motzkin and Robinson renewed interest in these
polynomials, and two more examples  in the style
of $M$ were  presented by Choi and Lam \cite{CL1,CL2}:
\begin{equation}
S(x,y,z) = x^4y^2 + y^4z^2 + z^4x^2 - 3x^2y^2z^2 \in \Delta_{3,6}, \\
\end{equation}
\begin{equation}
Q(x,y,z,w) = x^2y^2 + x^2z^2 + y^2z^2 + w^4 - 4wxyz \in \Delta_{4,4}.
\end{equation}

Here is an overview of the rest of the paper. 

In section two, we
present some preliminary material, mainly from curve theory; it is important
to consider reducible (as well as irreducible) polynomials. 

In section three, we present our version of Hilbert's Method (see
Theorem 3.4), based on more general perturbations and contradictions. 
There is a class of perturbations of a given positive polynomial with
fixed zeros by a polynomial which is singular at these zeros, in which
positivity is preserved. By counting dimensions, under certain
circumstances, there are  
polynomials of degree $2d$ which are singular on a set $A$, but are
not in the vector space generated by products of pairs of polynomials of degree
$d$ which vanish on $A$. If such a polynomial is positive, it cannot be sos. 
In Robinson's work, the set of cubics vanishing
at the eight points is spanned by $\{f_1,f_2\}$, but the vector space of
sextics which are singular at the eight points has dimension four and
so cannot be spanned by  $\{f_1^2,f_1f_2,f_2^2\}$. It is not necessary
to construct $\phi$ and $\psi$ to find this new sextic, although its behavior
at the ninth point must be analyzed to show that
 a successful perturbation is possible.

We show in Theorem 4.1 that Hilbert's Method works
when $f$ and $g$ are ternary cubics with exactly nine real
intersections, whether or not three are on a line or six on a
quadratic. (In other
words, Robinson's ``astonishing simplification'' always works.)
We also show that 
Hilbert's Method applies to the set of cubics which vanish on a
set of seven zeros, no four on a line, not all on a quadratic; see
Theorem 4.3.
\begin{example}
Let
\begin{equation}
\begin{gathered}
{\mathcal A} =
\{(1,0,0),(0,1,0),(0,0,1),(1,1,1),(1,1,-1),(1,-1,1),(1,-1,-1)\}, \\ 
F_1(x,y,z) = x(y^2-z^2), F_2(x,y,z) = y(z^2-x^2), F_3(x,y,z)
= z(x^2-y^2), \\ 
 G(x,y,z) = (x^2-y^2)(x^2-z^2)(y^2-z^2).
\end{gathered}
\end{equation}
It is easy to show that the $F_k$'s span the set of ternary cubics
which vanish on $\mathcal A$ and that $G$ is singular on $\mathcal A$
and not in the 
span of the $F_jF_k$'s. It follows from Theorem 4.3  that for some
$c>0$, $P_c = 
F_1^2+F_2^2+F_3^2+cG$ is psd and not sos. In fact, $P_1 = 2S$,
providing a new construction of (1.4).
\end{example}

In section five, we look at the sections of the cones  $P_{3,6}$ and
$\Sigma_{3,6}$ consisting of ternary sextics with the eight zeros of
Theorem 4.1. In addition to some general results,
we give a one-parameter family $\{R_t: t > 0\}$ of forms in
$\Delta_{3,6}$ with ten zeros and such that $R_1 = R$:
\begin{equation}
\begin{gathered}
R_t(x,y,z) :=  \\ 
\left(\frac{t^4+2t^2-3}{3} \right)(x^3-x z^2)^2 +
\left(\frac{1+2t^2-3t^4}{3t^4}\right)(y^3-y z^2)^2 + R(x,y,z).
\end{gathered}
\end{equation}
We give necessary and sufficient conditions for a sextic polynomial
$p(x,y)$ with zeros at $\{-1,0,1\}^2 \setminus (0,0)$ to be psd and to
be sos.

In section six, we present more examples in $\Delta_{3,6}$. This paper
would not be complete without an explicit illustration of Hilbert's
Method under his original restrictions. Theorems 4.1 and 4.3 and other
techniques  are then applied to produce new forms in
$\Delta_{3,6}$, including one-parameter families which
include $R$, $S$ and $M$.  For $t^2 < \frac 12$, let
\begin{equation}
\begin{gathered}
M_t(x,y,z) = (1-2t^2)(x^4y^2+x^2y^4) + t^4(x^4z^2+y^4z^2)\\- (3 -
8t^2+2t^4)x^2y^2z^2  -2t^2(x^2+y^2)z^4 + z^6;
\end{gathered}
\end{equation}
$M_t \in \Delta_{3,6}$ has ten zeros and $M_0 = M$.  Let
\begin{equation}
\begin{gathered}
S_t(x,y,z) = t^4(x^6+y^6+z^6) + (1-2t^6)(x^4y^2+y^4z^2+z^4x^2)\\ +
(t^8 - 2t^2)(x^2y^4+y^2z^4+z^2x^4)-3(1-2t^2+t^4-2t^6+t^8)x^2y^2z^2;
\end{gathered}
\end{equation}
$S_t \in \Delta_{3,6}$ has ten zeros if $t>0$. Note that $S_0 =
S$ and $S_1 = R$, so $S_t$ 
provides a ``homotopy'' between $S$ and $R$ in
$\Delta_{3,6}$ in the set of forms with ten zeros.  We also show that
\begin{equation}
U_c(x,y,z) = x^2y^2(x-y)^2 + y^2z^2(y-z)^2 + z^2x^2(z-x)^2 +
cxyz(x-y)(y-z)(z-x)
\end{equation}
is psd if and only if $|c| \le 4\sqrt{\sqrt 2 - 1}$ and sos only if $c=0$.
We conclude the section by returning to a subject brought up by Robinson:
$(ax^2+by^2+cz^2)R(x,y,z)$ is sos if and only if $a,b,c \ge
0$ and $\sqrt a, \sqrt b, \sqrt c$ are
 the sides of a (possibly degenerate) triangle.

 In section seven, we discuss the zeros of extremal ternary forms,
 using the perturbation argument from Hilbert's Method and show that
 if $p \in \Delta_{3,6}$ has exactly ten zeros, then it is
 extremal in the cone $P_{3,6}$. We present supporting evidence for the
 conjecture that, at least in a limiting sense, all extremal forms in
 $\Delta_{3,6}$ have ten zeros.

Finally, in section eight, we apply Hilbert's Method to provide a
family of positive polynomials in two variables in even degree $\ge 6$
which are not sos. We also 
speculate on the general applicability of Hilbert's Method in
higher degree.

Bezout's Theorem becomes more complicated in more variables,
and for that reason, we have confined our
discussions to ternary forms. However, we wish to record 
a somewhat unexpected connection between $\tilde R$ and $Q$
(c.f. (1.3), (1.5)):
\begin{equation}
\tilde R(x-w,y-w,z-w,x+y+z-w) = 2 Q(x,y,z,w).
\end{equation}
Robinson's example, after homogenization and this change in
variables, gives a new derivation of the Choi-Lam example. 
The set of quaternary quadratics which vanish on 
\begin{equation}
\begin{gathered}
{\mathcal A} =
\{(1,0,0,0),(0,1,0,0),(0,0,1,0),\\ (1,1,1,1),(1,1,-1,-1),
(1,-1,1,-1),(1,-1,-1,1)\}  
\end{gathered}  
\end{equation}
is spanned by $\{xy - zw, xz - yw, xw - yz\}$, and any such quadratic
also vanishes at $(0,0,0,1)$. The form $Q$ is evidently psd by the
arithmetic-geometric inequality,  singular on $\mathcal A$ and
positive at $(0,0,0,1)$, and so is not sos.

Parts of this paper have been presented at many conferences over the
last several years. The author thanks the organizers for their many
invitations to speak, and his friends and colleagues for
their encouragement and suggestions.

\section{Preliminaries}

Throughout this paper, we toggle between forms $F$ in $k$ variables
and polynomials $f$ in $k-1$ variables, with the ordinary convention that 
\begin{equation}
\begin{gathered}
f(x_1,\dots,x_{k-1}) := F(x_1,\dots,x_{k-1},1), \\
F(x_1,\dots,x_k) := x_k^d f(\tfrac{x_1}{x_k}, \dots,   \tfrac{x_{k-1}}{x_k}),
\end{gathered}
\end{equation}
where $d = \deg f$. For even $d$, it is easy to see that $F$ and $f$
are simultaneously psd or sos. It is usually more convenient to use
forms, since $F \in P_{k,m}$ if and only if $F(u) \ge 0$ for $u$ in
the compact set 
$S^{k-1}$, simplifying perturbation. On the other
hand, the zeros of $f$ can be isolated, whereas those of $F$ are not.

Following \cite{CLR1}, we define the {\it zero-set} of any $k$-ary $m$-ic
form $F$ by
\begin{equation}
\mz(F):= \{(a_1,\dots,a_k) \in \mathbb R^k\ : \ F(a_1,\dots,a_k) =0 \}.
\end{equation}
We have $0 \notin \mz(F)$ by convention, $|\mz(F)|$ will be
 interpreted as the number of lines in $\mz(F)$ and only
 only one representative of each line need be given. 
If $a \in \mz(F)$ and $a_k \neq 0$, then $a$ corresponds to
 a unique zero of $f$; if $a_k = 0$, then $a$
corresponds to a {\it zero of $f$ at infinity}. We also define
\begin{equation}
\mz(f):= \{(a_1,\dots,a_{k-1}) \in \mathbb R^{k-1}\ : \
f(a_1,\dots,a_{k-1}) =0 \}, 
\end{equation}
for non-homogeneous $f$. It is
possible for a strictly  positive $f$ to be have zeros at infinity. 
Consider $f(x,y) = x^2 + (xy-1)^2$ (and $F(x,y,z) = x^2z^2 + (xy-z^2)^2$):
clearly, $f(a,b)>0$ for $(a,b) \in \mathbb R^2$ and  $\mz(F) =
\{(1,0,0),(0,1,0)\}$.   

If $f$ is positive and $a \in \mz(f)$, then of course $\frac{\partial
  f}{\partial x_i}(a) = 0$  
for all $i$. We shall say that $f$ is {\it round at $a$} if $f_a$, the
second-order component of the
Taylor series to $f$ at $a$, is a positive definite quadratic
form. This is a ``singular non troppo''
 zero for a positive polynomial. The
corresponding second-order component of Taylor series for $F$ is psd
but not positive definite, since $F$ vanishes on lines through the
origin. 

If $F \in P_{n,m}$ (resp. $\Sigma_{n,m}$), and $G$ is derived from $F$
by an invertible linear change of variables, then  $G \in P_{n,m}$
(resp. $\Sigma_{n,m}$). Thus, it is harmless to assume when convenient
that $\mz(F)$ avoids the hyperplane $a_{n} = 0$; that is, $f$ has no zeros at
infinity.

Let $\mathbb R_{n,d} \subset \mathbb R[x_1,\dots,x_n]$ denote the
$\binom{n+d}n$-dimensional vector space of real polynomials
$f(x_1,\dots,x_n)$ with $\deg f \le d$. Suppose 
$A = \{\pi_1,\dots,\pi_r\} \subset \mathbb R^n$ is given. Let
$I_{s,d}(A)$ denote the vector space of those $p \in \mathbb R_{n,d}$ 
which have an 
$s$-th order zero at each $\pi_j$. In particular,
\begin{equation}
\begin{gathered}
I_{1,d}(A) = \{ p \in \mathbb R_{n,d}\ : \ p(\pi_j) = 0, \quad 1 \le j
\le r \}; \\
I_{2,2d}(A) = \left\{ p \in \mathbb R_{n,2d}\ : \ p(\pi_j) = \tfrac {\partial
  p}{\partial x_i}(\pi_j) = 0, \quad 1 \le i \le n,\quad 1 \le j \le r
\right\}. 
 \end{gathered}
\end{equation}
Since an $s$-th order zero in $n$ variables  imposes $\binom{n+s-1}n$
linear conditions, 
\begin{equation}
\dim I_{s,d}(A) \ge  \binom{n+d}n - r\binom{n+s-1}n.
\end{equation}
In Hilbert's sextic construction, $A = \{\pi_1,\dots,\pi_9\}$ is the
set of common zeros of $f_1(x,y)$ and $f_2(x,y)$,
 and  $\dim(I_{1,3}(A)) = 2 > \binom 52 - 9\binom 22$.

Let
\begin{equation}
I_{1,d}^2(A): = \left\{ \sum f_ig_i\ : \ f_i, g_i \in I_{1,d}(A) \right\}.
\end{equation}
Clearly, $I_{1,d}^2(A) \subseteq I_{2,2d}(A)$.
It is essential to Hilbert's Method that this inclusion may be strict;
for example, $\phi\psi(\pi_9) > 0$ so 
$\phi\psi \in  I_{2,6}(A) \smallsetminus I_{1,3}^2(A)$.

We also need to consider the ``forced'' zeros, familiar from the
Cayley-Bacharach Theorem; see \cite{EGH}.  Suppose 
$A\subset \mathbb R^n$ and $I_{1,d}(A)$ are given as above. Let 
\begin{equation}
\tilde A :=  \bigcap_{j=1}^r \mz(f_j) 
\smallsetminus A = \mz \biggl(\sum_{j=1}^r f_j^2 \biggr)
\smallsetminus A.  
\end{equation}

Unfortunately, this notation fails to capture forced zeros at
infinity. Accordingly, for $A \subset \mathbb R^n$, define the
associated projective set ${\mathcal A} \subset \mathbb R^{n+1}$ by
\begin{equation}
(a_1,\dots,a_n) \in A \iff (a_1,\dots, a_n,1) \in {\mathcal A}.
\end{equation}
As before, we define $I_{s,d}({\mathcal A})$ to be the set of $d$-ic
forms $F(x_1,\dots,x_{n+1})$ which have $s$-th order zeros on
$\mathcal A$. Then $f \in I_{s,d}(A)$ if and only if $F \in I_{s,d}({\mathcal
  A})$. We  define 
\begin{equation}
\tilde {\mathcal A} :=  \bigcap_{j=1}^r \mz(F_j) 
\smallsetminus {\mathcal A} = \mz \biggl(\sum_{j=1}^r F_j^2 \biggr)
\smallsetminus {\mathcal A}. 
\end{equation}
Given $A \subset \mathbb R^n$,  $\tilde{\mathcal A} = \emptyset$
when there are no forced zeros, even at infinity.

We say that $I_{1,d}(A)$ is {\it full} if, for any $\pi \in A$ and $v
\in \mathbb R^n$, there exists $f\in I_{1,d}(A)$ such that $\vec\nabla
f(\pi) = v$. Equivalently, if $\{f_1,\dots,f_s\}$ is a basis for
$I_{1,d}(A)$ and $f = \sum_j f_j^2$, then $I_{1,d}(A)$ is full if and
only if $f$ is round at each $\pi \in A$.

Bezout's Theorem in a relatively simple form is essential to our
proofs.  Suppose
$f_1(x,y)$ and $f_2(x,y)$ are relatively prime polynomials of degrees
$d_1$ and $d_2$. Let $\mz \subset \mathbb C^2$ denote the set of
common (complex) zeros of $f_1$ and $f_2$. Then
\begin{equation}
d_1d_2 = \sum_{\pi \in \mz} \mi_\pi(f_1,f_2),
\end{equation} 
where $\mi_\pi(f_1,f_2)$ measures the singularity of the intersection
of the curves $f_1=0$ and $f_2=0$ at $\pi$. 
In particular,  $\mi_\pi(f_1,f_2) = 1$ if and only if the curves
$f_1=0$ and $f_2=0$ are nonsingular at $\pi$ and have different tangents.
Thus, $\mi_\pi(f_1,f_2) = 1$ if and only if
$f_1^2+f_2^2$ is round at $\pi$, 
and $\mi_\pi(f_1,f_2) \ge 2$ otherwise. If $f_1$ and $f_2$ are both
singular at $\pi$, then $\mi_\pi(f_1,f_2) \ge 4$. 

\begin{lemma}
Suppose $f_1(x,y), f_2(x,y) \in \mathbb R_{2,d}$ and $|\mz(f_1) \cap
\mz (f_2)| = d^2$. If $A \subseteq \mz(f_1) \cap \mz (f_2)$ is
such that $I_{1,d}(A)$ has basis $\{f_1,f_2\}$, then $A$ is full.
\end{lemma}
\begin{proof}
It follows from (2.10) that any common zero of $f_1$ and $f_2$
must be real, and that $\mi_\pi(f_1,f_2) = 1$ for each
common zero $\pi$. It follows that $A$ is full.
\end{proof}

The next proposition collects some useful information from curve
theory. As is customary, if $f(\pi) = 0$, we say that  
 $\pi$ {\it lies on $f$}  or $f$ {\it contains} $\pi$.

\begin{proposition}
All polynomials herein 
are assumed to   be in $\mathbb R[x,y]$, and all 
enumerated sets of points are assumed to be distinct. These 
results apply to ternary forms with the obvious modifications.

\begin{enumerate}

\item If a quadratic $q$ is singular at $\pi$ and $q(\pi') = 0$ for
  some $\pi' \neq \pi$,   then $q = \ell_1\ell_2$ is a product of two
  linear forms $\ell_j$ containing $\pi$.

\item   If a set of eight points $A = \{\pi_1,\dots,\pi_8\}$
  is given, no four   on a line and no seven on a quadratic, then
  $\dim I_{1,3}(A) = 2$. 

\item In the last situation, if $A_j = A \smallsetminus \{\pi_j\}$, then
  there exists a cubic $f$ so that $f |_{A_j} = 0$, but $f(\pi_j) \neq
  0$; in particular,  $\dim I_{1,3}(A_j) = 3$.

\item Suppose $f(x,y)$ and $g(x,y)$ are cubics, $A = \mz(f) \cap
  \mz (g) = \{\pi_1,\dots,\pi_9\}$ and $A_j = A \smallsetminus  \{\pi_j\}$. 
 For each $j$, 
  $I_{1,3}(A_j) = I_{1,3}(A)$. In other words, if eight of the points
  lie on a cubic $h$, then so will the ninth.

\item Under the same conditions as (4), no four of the $\pi_i$'s
  lie on a line and no seven lie on a quadratic. Three of the
  $\pi_i$'s lie on a line if and 
  only if the other six  lie on a quadratic if
  and only if  $I_{1,3}(A)$ contains a reducible cubic.

\end{enumerate}
\end{proposition}
\begin{proof}
For (1), write $q(x,y) = a + bx + cy + dx^2 + exy + fy^2$ and assume by
translation that $\pi = (0,0)$. Then $a=b=c=0$ and $q(x,y) = dx^2 +
exy + fy^2$. If $\pi' = (r,s) \neq (0,0)$, then $sx - ry$ is a factor
of $q$. The next two assertions are classical and proofs can be found,
for example, in \cite[Ch.15]{Bi}; (4) is well-known and  is often
attributed to Cayley-Bacharach,  but it was discovered by Chasles; see \cite{EGH}.

 For (5),
if four $\pi_i$'s lie on a line $\ell$, then $\ell$ divides both
$f$ and $g$ by Bezout, so that $|\mz(f) \cap \mz (g)| = \infty$. If seven
$\pi_i$'s lie on a reducible quadratic $q = \ell_1\ell_2$, then at
least four lie on one  $\ell_i$, and we are in the earlier case. If
they lie on an irreducible $q$, 
then it must be indefinite, and again, $q$ divides both
$f$ and $g$ by Bezout, so that $|\mz(f) \cap \mz (g)| = \infty$. 

Suppose now that three points of $A$, say $\{\pi_1,\pi_2,\pi_3\}$, lie
on the line $\ell$ and let   
$q$ be the quadratic containing  $\{\pi_4,\dots,\pi_8\}$. Then $\ell q =
0$ on $A_8$, so by (4), $(\ell q)(\pi_9) = 0$. Since $\ell(\pi_9)\neq 0$,
we must have $q(\pi_9) = 0$; thus six zeros lie on $q$ and $\ell q \in
I_{1,3}(A_j)$. 
(A similar proof follows if we start with six points lying on the
quadratic $q$.) Finally, if $\ell q \in
I_{1,3}(A)$, then at most three of the $\pi_i$'s can lie on $\ell$, and
at most six can lie on $q$, hence these numbers are exact.
\end{proof}

\begin{lemma}
Suppose $A$ is a set of eight distinct points, no four on a line and no seven
on a quadratic, and let $\{f_1,f_2\}$ be a basis for
$I_{1,3}(A)$. Then $f_1$ and $f_2$ are relatively prime.
\end{lemma}
\begin{proof}
 If $f_1$ and $f_2$ have a common quadratic factor $q$,
then $f_i = \ell_i q$ and at most six points of $A$ lie on $q$, so
$\ell_1$ and $\ell_2$ share two points and so are proportional, a
contradiction. If $f_1$ and $f_2$ have only  a common linear factor $\ell$,
then $f_i = \ell q_i$, and at most three points of $A$ lie on $\ell$, so
$q_1$ and $q_2$ share five points and so are proportional, again a
contradiction. 
\end{proof}

In the situation of Lemma 2.3, 
Bezout's Theorem has one of three possible implications:
(a) there is a ninth point $\pi \in \tilde A$ so that $f_1(\pi) = f_2(\pi) =
0$; (b) $\tilde A = \emptyset$, but $(a,b,0) \in \tilde {\mathcal A}$
is a common zero of $f_1$ and $f_2$ at infinity; (c)
$\mi_\pi(f_1,f_2) = 2$ 
for some $\pi \in A$. The first two cases are essentially the same: if (b)
occurs, we homogenize and change variables so that the zero is no
longer at infinity after dehomogenization. Any necessary construction
can then be performed, and the variables changed back.
The third case is singular, but seems to be difficult to identify in
advance, and is equivalent to the existence of a cubic in $I_{1,3}(A)$
which is singular at some $\pi \in A$.

We shall say that a set of eight points $A$ for which (a) or (b)
occurs is {\it  
  copacetic}. Since $f_1$ and $f_2$ are real, $f_1(\pi) = f_2(\pi) = 0
 \implies  f_1(\bar \pi) = f_2(\bar \pi) = 0$. Bezout implies that
 $\pi = \bar \pi$; 
that is, the ninth point $\pi$ must be real. We have the following corollary
to Lemma 2.1.

\begin{lemma}
If $A$ is copacetic, then it is full.
\end{lemma}

The following lemma was probably known a hundred years ago.

\begin{lemma} 
Suppose seven points $A = \{\pi_1,\dots,\pi_7\}$ in the plane are given,
not all on a quadratic and no four on a line. Then, up to multiple,
 there is a unique cubic $f(x,y)$ which is singular at $\pi_1$ and contains
 $\{\pi_2,\dots,\pi_7\}$.
\end{lemma}
\begin{proof}
Since $1 \cdot 3+6\cdot 1< 10$ linear conditions are given, at least one such
nonzero $f$
exists. Suppose $f_1$ and $f_2$  satisfy these properties and are not
proportional.  Then 
$\sum_j \mi_{\pi_j}(f_1,f_2) \ge 2^2 + 6\cdot 1 > 3 \cdot 3$, hence
$f_1$ and $f_2$ have a common factor. The common factor could be an
irreducible quadratic, a reducible quadratic, or linear.

In the first case, $f_1 = \ell_1 q$ and $f_2 = \ell_2 q$, where
$q(\pi_1) = \ell_i(\pi_1) = 0$ by Prop. 2.2(1). At least one point,
say $\pi_7$, does not lie on $q$, hence $\ell_i(\pi_7) = 0$ as well.
Thus the two $\ell_i$'s share two zeros and are proportional, a contradiction.

In the second case, we have $f_1 = \ell_1\ell_2\ell_3$ and $f_2 =
\ell_1\ell_2\ell_4$, and $\pi_1$ lies on at least two of 
$\{\ell_1,\ell_2,\ell_3\}$ and two of
$\{\ell_1,\ell_2,\ell_4\}$. If $\ell_1(\pi_1) = \ell_2(\pi_1) = 0$,
then $\ell_1$ and $\ell_2$ together can contain at most four of the six points
$\{\pi_2,\dots,\pi_7\}$, hence $\ell_3$ and $\ell_4$ must each contain
at least two points in common, and so are proportional, again
a contradiction. Otherwise, without loss of generality,
$\ell_1(\pi_1)=0$ and $\ell_2(\pi_1) \neq 0$, hence
$\ell_3(\pi_1)=\ell_4(\pi_1)=0$. 
In this case,  $\ell_1$ and $\ell_2$ can
together contain at most five of the six points
$\{\pi_2,\dots,\pi_7\}$, so that 
$\ell_3$ and $\ell_4$ must contain also some $\pi_j$ other than $\pi_1$. This
is again a contradiction.

Finally, suppose $f_1 = \ell q_1$ and $f_2 = \ell q_2$, where $q_1$
and $q_2$ are relatively prime quadratics, so they share at most four
points.  If $\ell(\pi_1) = 0$, then $q_j(\pi_1) = 0$ as well and since
at least four of   $\{\pi_2,\dots,\pi_7\}$ do not lie on $\ell$, they
must lie on both $q_1$ and $q_2$. Thus $q_1$
and $q_2$ share five points, a contradiction. If $\ell(\pi_1)
\neq 0$, then $h_1 = \ell \ell_1 \ell_2$ and $h_2 = \ell \ell_3 \ell_4$,
where the $\ell_i$'s are distinct lines containing $\pi_1$. But if
$\ell(\pi_j) \neq 0$ (which is true for at least four $\pi_j$'s) then
$\pi_j$ must also lie on one of $\{\ell_1,\ell_2\}$ and one of
$\{\ell_3,\ell_4\}$. That is, the line through $\pi_1$ and $\pi_j$
divides both $\ell_1\ell_2$ and $\ell_3\ell_4$, a final contradiction.
\end{proof}

The last lemma in this section is used in the proof of Theorem 4.3. 

\begin{lemma}
If d=3 and $A$ is a set of seven points in $\mathbb R^2$, no four on a line
and not all on a quadratic, then $A$ is full and
$\tilde{\mathcal A} = \emptyset.$
\end{lemma}
\begin{proof}
Choose $\pi_8$ to avoid any line between two points of $A$ and any
quadratic determined by five points of $A$. Then $A \cup \{\pi_8\}$
has no four points in a line and no seven on a quadratic, and so 
 $\dim I_{1,3}(A) = 3$ by Prop. 2.2(3). Suppose $\{f_1,f_2,f_3\}$ is
a basis for $I_{1,3}(A)$ and for each $j$, consider the map
\begin{equation}
T_j: (c_1,c_2,c_3) \mapsto \sum_{k=1}^3 c_k \vec\nabla f_k(\pi_j).
\end{equation}
By Lemma 2.5, $\dim(\ker(T_j)) = 1$, hence each $T_j$ is surjective, and 
so $A$ is full.

Suppose $\pi \in \tilde {\mathcal A}$; after an invertible linear
change, we may assume without loss of generality that $\pi \in \tilde
A$. By the contrapositive to Prop. 2.2(3), either $A \cup
\{\pi\}$ has four points in 
a line or has seven points on a quadratic.
Again, choose $\pi_8$ so that $A_1 = A \cup \{\pi_8\}$ has no four
points in a line and no seven on a quadratic. By Prop. 2.2(2),
we may assume without loss of generality that $I_{1,3}(A_1)$ has basis
$\{f_1,f_2\}$, so  $\pi \in \tilde A_1$. Let $A_2 = A_1 \cup
\{\pi\}$.  Thus $f_1$ and $f_2$ are two cubics which vanish on a set
$A_2$ with four points on a line $\ell$ or seven points on a quadratic
$q$, and so $f_1$ and $f_2$ have a common factor by Bezout, a
contradiction.
\end{proof}

\section{Hilbert's Method}

We begin this section with  a general perturbation result.

\begin{lemma}[The Perturbation Lemma]
Suppose $f,g \in \mathbb R_{n,2d}$ satisfy the following conditions:

\begin{enumerate}
\item The polynomial $f$ is positive with no zeros at infinity,
  and $2d = \deg f \ge \deg g$; 

\item There is a finite set $V_1$ so that if $v \in
  V_1$, then $f$ is round at $v$ and $g$ vanishes to second-order at $v$; 

\item The set $V_2:= \mz(f)\smallsetminus V_1$ is finite and if $w \in
  V_2$, then $g(w) > 0$.

\end{enumerate}
Then there exists $c = c(f,g)>0$ so that $f+cg $ is a positive polynomial. 
\end{lemma}
\begin{proof}
For $v \in V_1$, let $g_v$ denote the second-order (lowest degree)
term of the Taylor
series for $g$ at $v$. Since $f_v$ is positive definite, there 
exists $\al(v) > 0$ so that $f_v + \al g_v$ is positive definite for $0
\le \al \le \al(v)$. 
If $\al_0 = \min_v \al(v)$, then there exist neighborhoods
$\mathcal N_v$ of each $v$  so that $f + \al_0 g$ is
positive on each ${\mathcal N_v} \smallsetminus \{v\}$.
 Further, for $w\in V_2$, $(f+\al_0 g)(w) = \al_0 g(w) > 0$, 
hence there is a neighborhood $\mathcal N_w$ of $w$ on which  $f + \al_0
g$ is positive. It follows that $f + \al_0g$ is non-negative on the open set 
$\mathcal N = \cup \mathcal N_v \cup \mathcal N_w$. 

Homogenize $f,g$ to forms $F,G$ of degree $2d$ in $n+1$. For
$x \in \mathbb R^{n}$, let $||x||
= (1+\sum_i x_i^2)^{1/2}$ and let $\widetilde {\mathcal N}$ be 
 the image of $\mathcal N$ under the map
\begin{equation}
(x_1,\dots,x_n) \mapsto \left( \frac{x_1}{||x||},\dots,
  \frac{x_n}{||x||}, \frac 1{||x||}\right)  \in S^n. 
\end{equation}
Then $\widetilde {\mathcal N}$ is open and $(F+\al G)(x) \ge 0$ for $x
\in \widetilde {\mathcal N}$. By hypothesis, $\mz(F) \subset
\widetilde {\mathcal N}$, hence $F$ is positive on the complement
$\widetilde{\mathcal N}^c$, so it achieves a positive minimum on the
compact set 
$\widetilde{\mathcal N}^c$. Since $G$ is bounded on $S^n$, there exists
$\beta > 0$ so that $(F+\beta G)(x) \ge 0$ for $x \in \widetilde
{\mathcal N}^c$. It follows that
$F+cG$ is psd, where $c = \min\{\al_0,\be\}$. The desired result follows
upon dehomogenizing. 
\end{proof}

The following two theorems generalize the contradiction of Hilbert's
construction.

 \begin{theorem}
If $p \in I_{2,2d}(A)$ is sos, then $p \in I_{1,d}^2(A)$.
\end{theorem}
\begin{proof}
If $p = \sum_k h_k^2$, then $p(a) = 0$ for $a \in A$, hence  $h_k(a) =
0$, and so  $h_k \in I_{1,d}(A)$, implying $p \in I_{1,d}^2 (A)$.
\end{proof}  

Let $I_{1,d}(A)$ have basis $\{f_1,\dots,f_r\}$, and suppose 
the $\binom {r+1}2$ polynomials $f_if_j, 1 \le i \le j \le r$ are 
linearly independent; in other words, for each $p \in I_{1,d}^2(A)$ there is
a unique quadratic form $Q$ so that $p = Q(f_1,\dots,f_r)$. We call
this the {\it independent case}. (We have been unable to find $I_{1,d}(A)$
for which this does not hold.)  Let 
\begin{equation}
R_f:= \{ (f_1(x),\dots,f_r(x))\ : \ x \in \mathbb R^n \} \subseteq \mathbb
R^r.
\end{equation}
denote the range of the basis as an $r$-tuple. 

\begin{theorem}
Suppose  $p = Q(f_1,\dots,f_r) \in I_{1,d}^2(A)$ in the independent case:
\begin{enumerate}
\item $p$ is sos if and only if $Q$ is an sos quadratic form;
\item $p$ is psd if and only if $Q(u) \ge 0$ for $u \in R_f$;
\item if $n=2$, $r=2$, and $f_1$ and $f_2$ are relatively prime
  polynomials with odd  degree $d$, then $R_f
  = \rr^2$, hence $p \in   I_{1,d}^2(A)$ is psd if and only if it is sos. 
\end{enumerate}
\end{theorem}
\begin{proof}
If $p = \sum_k h_k^2$ is sos, then as in the last proof, $h_k \in
I_{1,d}(A)$. To be specific, if $h_k = \sum_\ell
c_{k\ell} f_\ell$, then by the uniqueness of $Q$,
 $Q(u_1,\dots,u_r) = \sum_\ell
(\sum_\ell c_{k\ell}u_\ell)^2$. Conversely, if $Q = \sum_\ell T_\ell^2$ for
linear forms $T_\ell$, then $p = \sum_\ell T_\ell(f_1,\dots,f_r)^2$.

The assertion in (2) is immediate.

For (3), we first note that, since
$(f_1(\la x), f_2(\la x)) = \la^d (f_1(x),f_2(x))$, it suffices to
show that every line through the origin intersects $R_f$. By
hypothesis, $\mz(f_1)$ and $\mz(f_2)$ are infinite sets, but
$|\mz(f_1) \cap \mz(f_2)| \le d^2$.
It follows that there exist $\pi$ and $\pi'$ so that
$(f_1(\pi),f_2(\pi)) = (1,0)$ and $(f_1(\pi'),f_2(\pi')) = (0,1)$. Now
take a curve $\ga(t) \in \rr^2$ and so that $\ga(0) = \pi$, $\ga(1) =
\pi'$ and  $\ga(2) = -\pi$ and $\ga(t) \notin \mz(f_1) \cap
\mz(f_2)$, and let $h(t) =(f_1(\ga(t)),f_2(\ga(t)))$. We have $h(0) =
(1,0)$, $h(1) = (0,1)$, $h(2) = (-1,0)$ and $h(t) \neq (0,0)$, so by
continuity, 
each line through the origin contains some $h(t)$, $0 \le t \le 2$.
\end{proof}

The hypotheses of Theorem 3.3(3) applies in Hilbert's
original construction, with $d=3$. We show in Example 3.1 below
that $R_f \neq \rr^r$ in general, and combine this with Theorem 3.3(2)
to give one instance of a positive form in a $I_{1,d}^2(A)$  which is
not sos. 
\begin{theorem}[Hilbert's Method]
Suppose a finite set $A \subset \mathbb R^n$ is such that $I_{1,d}(A)$
has basis $\{f_1,\dots,f_s\}$, where $\tilde A$ is
finite,  $A$ is full and $f = \sum_j f_j^2$ has no zeros at infinity.
Further, suppose there exists $g \in I_{2,2d}(A) \smallsetminus
I_{1,d}^2(A)$ so that  $g(w) > 0$ for each $w \in \tilde A$. Then there
exists $c > 0$ so that 
\begin{equation}
p_c = \sum_{j=1}^s f_j^2 + c g
\end{equation}
is positive and not sos.
\end{theorem}
\begin{proof}
In the notation of Lemma 3.1, let $V_1 = A$ and $V_2 = \tilde A$.  
Since $f$ has no zeros at infinity, $\deg f = 2d$, and 
$A$ is full, the hypotheses of Lemma 3.1 are satisfied.
Thus there exists $c>0$ so that $p_c$ is positive, and since $p_c \notin
I_{1,d}^2(A)$, it is not sos by Theorem 3.2.
\end{proof}

\begin{remarks}
\qquad

\begin{enumerate}

\item
 If $\tilde A = \emptyset$,
then the Perturbation Lemma can be 
applied to $(f,\pm g)$ for both signs, so that $f \pm c g$ is positive
for some $c > 
0$ and  both choices of sign.

\item
In any particular case, the condition that $f$ is round at $v \in V_1$
can be relaxed in the Perturbation Lemma, so long as a stronger
condition is imposed  on $g$ to insure  
that $f + \al g$ is positive in some punctured 
neighborhood $\mathcal N_v$ of $v$.

\item
Since Hilbert's Method applies to any basis of $I_{1,d}(A)$, we may
replace $\sum_j f_j^2$ by any positive definite quadratic form in the
$f_j$'s. 

\item
Hilbert's original sextic contradiction follows from
$(\phi\psi)(\pi_9) \neq 0$, which implies that $\phi\psi \in
I_{2,6}(A)\smallsetminus I_{1,3}^2(A)$.

\item
 Theorem 4.3 covers a situation
in which  $\tilde A = \emptyset$, but that $I_{2,2d}(A)\smallsetminus
I_{1,d}^2(A)$ is non-empty, so Hilbert's Method still applies.
\end{enumerate}
\end{remarks}

\begin{example}
We revisit Example 1.1, keeping the notation of (1.6). It is easy to
check that  
\begin{equation}
\{F_1^2,F_2^2,F_3^2,F_1F_2,F_1F_3,F_2F_3\}
\end{equation}
is linearly independent, so that Theorem 3.3 applies.
Let
\begin{equation}
Q(u_1,u_2,u_3) = 5u_1^2 + 5u_2^2 + 5u_3^2 -6u_1u_2-6u_1u_3-6u_2u_3;
\end{equation}
evidently, $Q$ is not a psd quadratic form. We show now (in two ways) that
\begin{equation}
T:=  Q(F_1,F_2,F_3)
\end{equation}
is psd; note that $T$ is not sos by Theorem 3.3(1).

Let
\begin{equation}
\begin{gathered}
P(v_1,v_2,v_3) : =  v_1^4 +v_2^4+v_3^4 -
2v_1^2v_2^2-2v_1^2v_3^2-2v_2^2v_3^2 \\ =
(v_1+v_2+v_3)(v_1+v_2-v_3)(v_1-v_2+v_3)(v_1-v_2-v_3).
\end{gathered}
\end{equation}
A computation shows that
\begin{equation}
P(F_1,F_2,F_3) = (x^2-y^2)^2(x^2 - z^2)^2(y^2 - z^2)^2
\end{equation} 
is psd, hence $R_F \subseteq \{(x,y,z): P(x,y,z) \ge 0\}$. We claim that
that $Q \ge 0$ on $R_F$ and so $T$ is psd by Theorem 3.3(2).
Since
\begin{equation}
 5u_1^2 + 5u_2^2 + 5u_3^2 -6u_1u_2-6u_1u_3
\end{equation}
is psd, if $\bar u_2\bar u_3 < 0$, say, then $Q(\bar u_1,\bar u_2,\bar
u_3) \ge 0$. By symmetry, it follows that
$Q(v_1,v_2,v_3)\ge 0$ unless the $v_i$'s have the same sign, and it
suffices to suppose $v_1\ge v_2 \ge v_3 \ge 0$. The first
three linear  factors of $P$ in (3.7) are always positive, so
$P(v_1,v_2,v_3) \ge 
0$ if and only if $v_1 = v_2 + v_3 + t$ with $t \ge 0$. Since 
\begin{equation}
Q(v_2+v_3+t,v_2,v_3) = 4(v_2-v_3)^2 + t(4v_2+4v_3+5t),
\end{equation}
the claim is verified.

The second proof is direct. We note that  $T$ is symmetric:
\begin{equation}
T(x,y,z) = 5\sum^6 x^4y^2 + 6\sum^3 x^4yz + 6\sum^3 x^3y^3 -
 6\sum^6 x^3y^2z - 30x^2y^2z^2.
\end{equation}
A calculation shows that
\begin{equation}
\begin{gathered}
2(x^2+y^2+z^2-xy-xz-yz)T(x,y,z) = 
(x-y)^4(xy + 3xz+3yz+z^2)^2 \\ + (x-z)^4(xz + 3xy+3yz+y^2)^2 +
(y-z)^4(yz + 3xy+3xz+x^2)^2,
\end{gathered}
\end{equation}
so $T$ is psd. Although $|\mz(T)| = 7$, the zeros at
$(1,1,-1),(1,-1,1),(-1,1,1)$ are not round. In fact, $T(1+t,1-t,-1) =
48t^4 + 4t^6$, etc. These singularities are
 useful in constructing the representation (3.12).
\end{example}

\section{Two applications of Hilbert's Method to ternary sextics}

In this section we show that Robinson's simplification of Hilbert's
Method works in general. By Theorem 2.2(5), the assumption that no
three of the nine points are on a line and no six are on a quadratic
is equivalent to saying that no $\al f_1 + \be f_2$ is reducible.
Theorem 4.1 removes this restriction. In Theorem 4.3, we show that
Hilbert's Method also applies to the set of ternary
sextics which share seven zeros, no four in a line, no seven on a
quadratic. 
\begin{theorem}
Suppose $f_1(x,y)$ and $f_2(x,y)$ are two
relatively prime real cubics with exactly nine distinct real
common zeros. Then Hilbert's Method applies to any subset $A$ of eight of
the common zeros. 
\end{theorem}
\begin{proof}
Lemma 2.4 shows that if $A = \{\pi_1,\dots,\pi_8\}$
is copacetic, as is assumed here, then $\tilde A =\{\pi_9\}$ and $A$ is full.
It follows from (2.5) that $\dim I_{2,6}(A) \ge \binom 82 - 3
\cdot 8 = 4$. Since $I_{1,3}^2(A)$ is spanned by $\{f_1^2,f_1f_2,f_2^2\}$,
there exists $0 \neq g \in I_{2,6}(A) \smallsetminus
I_{1,3}^2(A)$. If we can show that $g(\pi_9)
\neq 0$, then $\pm g(\pi_9) > 0$ for some choice of sign, and 
Theorem 3.4 applies. 

Suppose to the contrary that $g(\pi_9) = 0$. Either $g$ is
singular at $\pi_9$, or there exists $(\al_1,\al_2) \neq (0,0)$ so that
the tangents of $g$ and $\al_1f_1+\al_2f_2$ are parallel at
$\pi_9$. Since the choice of basis for $I_{1,3}(A)$ was arbitrary, we
may assume without loss of generality that $(\al_1,\al_2) = (1,0)$
from the beginning. In either case, $ \mi_{\pi_9}(f_1,g) \ge 2$, so
\begin{equation}
 \sum_{j=1}^9 \mi_{\pi_j}(f_1,g) \ge 2 \cdot 9 = \deg(f_1) \cdot \deg(g).
\end{equation}
Since $f_1$ is a real cubic, there exists $\pi_0 \notin A \cup \tilde
A$ so that $f_1(\pi_0) = 0$ 
and, necessarily, $f_2(g_0) \neq 0$. Now let
\begin{equation}
\tilde g = g - \frac{g(\pi_0)}{f_2^2(\pi_0)} f_2^2,
\end{equation} 
so that $\tilde g(\pi_0) = 0$.
Observe that  $\tilde g \in I_{2,6}(A) \smallsetminus I_{1,3}^2(A)$, and $g$
and $\tilde g$ agree to second-order at $\pi_9$. In particular, they
are either both singular or have the same tangents. Thus, we may
replace $g$ by $\tilde g$ for purposes of the argument, and
assume that $g(\pi_0) = 0$. Combining $\mi_{\pi_0}(f_1,g) \ge 1$ with
(4.1), we see that $f_1$ and $g$ have a common factor by Bezout. Let $d =
\deg(\gcd(f_1,g))$.

If $d=3$, then $g = f_1 k$ for some cubic $k$. Since $g$ is singular
on $A$ and $f_1$ is singular at no point of $A$, we must have $k \in
I_{1,3}(A)$, so that $g \in I_{1,3}^2(A)$, a contradiction. (Under
Hilbert's restrictions, $f_1$ is irreducible, so this is the
only case.)

Suppose $d = 2$ and write $f_1 = \ell q$ and $g = p q$, where $\ell$
is linear, $q$ is quadratic and $p$ is quartic and $\ell$ and $p$ are
relatively prime. Then $\ell=0$
on exactly three of the $\pi_i$'s. After reindexing, there are two
cases: either $\ell = 0$ on $\{\pi_1,\pi_2,\pi_3\}$ or $\ell = 0$ on
$\{\pi_1,\pi_2,\pi_9\}$, with $q = 0$ on the complementary sets.
In the first case, $q(\pi_i)\neq 0$ for $i=1,2,3$, so $p$ is singular
at these three points and $\mi_{\pi_1}(\ell,p) + \mi_{\pi_2}(\ell,p) +
\mi_{\pi_3}(\ell,p) \ge 6 > 1 \cdot 4$. Since $\ell$ and $p$ are
relatively prime, this is a contradiction by Bezout. In the second
case, $p$ is still singular at $\pi_1,\pi_2$ and $q(\pi_9) \neq 0$, so
$p(\pi_9) = 0$ and 
$\mi_{\pi_1}(\ell,p) 
+ \mi_{\pi_2}(\ell,p) +  \mi_{\pi_9}(\ell,p) \ge 5 > 2+2+1$,
another contradiction.

Finally, suppose $d = 1$ and write $f_1 = \ell q$ and $g = \ell p$,
where $\ell$ is linear, $q$ is quadratic and $p$ is quintic and $q$
and $p$ are relatively prime. With either case for $\ell$ as above,
$\ell\neq0$ and $p$ is singular at $\pi_4,\dots,\pi_8$ and $\mi_{\pi_4}(q,p) +
\dots + \mi_{\pi_8}(q,p) \ge 10 = 2 \cdot 5$. In the first case,
$\ell(\pi_9) \neq 0$, 
so $\mi_{\pi_9}(q,p) \ge 1$; in the second case,  $\ell(\pi_3) \neq 0$, so
$\mi_{\pi_3}(q,p) \ge 2$. In either case Bezout implies that $q$
and $p$ are not relatively prime, and this contradiction completes the
proof. 
\end{proof}

It is possible for $g$ and the $f_i$'s to have a common factor,
provided it does not contain $\pi_9$.
This happens in Robinson's example: $f_1 = x(x^2-1)$, $f_2
= y(y^2-1)$ and $g = (x^2-1)(y^2-1)(1-x^2-y^2)$.

\begin{corollary}
If $A$ is copacetic, then there exists a positive sextic polynomial
$p(x,y)$ so that $A \subseteq \mz(p)$ and $p$  is not sos.
\end{corollary}

\begin{theorem}
Suppose $A = \{\pi_1,\dots,\pi_7\}\subset \mathbb R^2$, 
with no four $\pi_i$'s in a line and
not all seven on one quadratic. Then Hilbert's Method applies to $A$.
\end{theorem}
\begin{proof}
It follows from Lemma 2.6 that
$A$ is full and $\tilde {\mathcal A} = \emptyset$. We have
$\dim I_{1,3}(A) = 3$, so $\dim I_{1,3}^2(A) \le 6$, but by (2.5),
$\dim I_{2,6}(A) \ge \binom 82 - 7\cdot \binom 32 = 7$. Thus there
exists $g \in \dim I_{2,6}(A) \smallsetminus I_{1,3}^2(A)$ and since
 $\tilde {\mathcal A} = \emptyset$, Hilbert's Method can be applied.
\end{proof}
Theorem 4.3 is implemented in Examples 1.1 and 6.3.
\begin{corollary}
If $A$ is a set of seven points in $\mathbb R^2$, no four on a line and not all
 on a quadratic, then there exists a positive sextic polynomial
$p(x,y)$ so that $A \subseteq \mz(p)$ and $p$  is not sos.
\end{corollary} 

\section{Psd and sos sections}

We now consider $I_{2,6}({\mathcal A}) \cap P_{3,6}$ and
$I_{2,6}({\mathcal A}) \cap \Sigma_{3,6}$ in detail.
Our motivation is that $P_{3,6}$ and $\Sigma_{3,6}$ lie
in $\mathbb R^{28}$ and are difficult to visualize. These two
sections, in general, lie in $\mathbb R^4$, and thus are
more comprehensible. We work in the homogeneous case.

\begin{theorem}
In the notation of Theorem 4.1, suppose
\begin{equation}
P = c_1 F_1^2 + 2 c_2 F_1F_2 + c_3 F_2^2 + c_4 G.
\end{equation}
If $P$ is sos, then $c_4 = 0$. If $c_4=0$, then $P$ is sos if and only
if $P$ is psd if and only if $c_1 \ge 0$, $c_3 \ge 0$ and $c_1c_3 \ge c_2^2$.
\end{theorem}
\begin{proof}
These are Theorems 3.2 and 3.3(1),(2) in the homogeneous case. 
\end{proof}

Because $G$ is only defined modulo $I_{1,d}^2({\mathcal A})$, it is
difficult to make any general statements about the circumstances under
which $P$ is psd.  However, one can identify the possible zeros of
$P$.
\begin{theorem}
Suppose $P = c_1 F_1^2 + 2 c_2 F_1F_2 + c_3 F_2^2 + c_4 G$ is psd,
where $c_4 \neq 0$ and let $J$ be the Jacobian of $F_1,F_2$ and $G$. Then
\begin{equation}
 \mz(P) \subseteq \mz(F_1) \cup \mz(F_2) \cup  \mz(J).
\end{equation}
\end{theorem} 
\begin{proof}
If $P(a) = 0$ and  $(F_1(a),F_2(a)) \neq (0,0)$, then $P$ and
$(F_1(a)F_2-F_2(a)F_1)^2$ are linearly independent sextics which are both
singular at $a$. Thus the Jacobian of $(F_1^2, F_1F_2, F_2^2,G)$, when
evaluated at $a$, has  rank $\le 2$. In
particular, the $3 \times 3$ minor omitting $F_1F_2$ vanishes; this
minor reduces to $4F_1F_2J$.
\end{proof}

A maximal perturbation
might not lead to a new zero, but rather to a greater singularity
at a pre-existing zero; see Example 3.1.

In the special case of Robinson's example, we are able to give a much
more precise description of these sections. Let $A =
\{-1,0,1\}^2 \smallsetminus \{(0,0)\}$. A routine calculation shows that
$f_1(x,y) = x^3-x$ and $f_2(x,y) = y^3-y$ span $I_{1,3}(A)$ and
$f_1^2,f_1f_2,f_2^2$ and $g(x,y) = 
(x^2-1)(y^2-1)(1-x^2-y^2)$ span $I_{2,6}(A)$. It is convenient to
replace $g$ with $f_1^2+f_2^2+g$, which homogenizes to $R$.

Consider now
\begin{equation}
\begin{gathered}
\Phi[c_1,c_2,c_3,c_4](x,y,z):= c_1F_1^2 + 2c_2F_1F_2+c_3F_2^2+c_4R
\\ = c_1(x^3-xz^2)^2 +
2c_2(x^3-xz^2)(y^3-yz^2)+c_3(y^3-yz^2)^2 \\
+ c_4(x^6 + y^6 + z^6 - x^4y^2 - x^2y^4 - x^4z^2-y^4z^2 -
x^2z^4-y^2z^4+3x^2y^2z^2).
\end{gathered}
\end{equation}
This is the general form of $\Phi \in I_{2,6}({\mathcal A})$, where
\begin{equation}
{\mathcal A}  = \{(\pm 1,0,1), (0,\pm 1,1),(\pm 1, \pm 1,1)\}.
\end{equation}
Theorem 5.1 implies that
$\Phi(c_1,c_2,c_3,0)$ is psd if and only if it is sos if and only if
$c_1,c_3,c_1c_3-c_2^2 \ge 0$, so we may henceforth assume that $c_4 \neq 0$.

We begin our discussion of positivity with a collection of short
observations. 

\begin{lemma}
Suppose $\Phi[c_1,c_2,c_3,c_4]$ is psd. Then the following are true:
\begin{enumerate}
\item $c_4 \ge 0$;

\item   $\Phi[c_1,-c_2,c_3,c_4]$ and  $\Phi[c_3,c_2,c_1,c_4]$ are psd;

\item $\Gamma(x,y):=
 (c_1+c_4) x^6 - c_4 x^4y^2 + 2c_2 x^3y^3 - c_4 x^2y^4 + (c_3+c_4)
  y^6$ is psd; 
 
\item $\Phi[c_1,0,c_3,c_4]$ is psd.
\end{enumerate}
\end{lemma}

\begin{proof}
The first observation follows from evaluation at $(0,0,1)$, the second
from taking $(x,y,z)\mapsto (x,-y,z),(y,x,z)$, the third from setting
$z=0$, and the fourth from averaging the psd forms $\Phi[c_1,\pm
  c_2,c_3,c_4]$. 
\end{proof}

In view of Lemma 5.3(1), it suffices now to assume $c_4 = 1$.
For $t > 0$, let
\begin{equation}
\al(t) = \frac{2t^2+t^4}3, \qquad \be(t) = \frac{1+2t^2}{3t^4}, \qquad
\ga(t) = \be(\al^{-1}(t)).
\end{equation}
Then $\be(t) = \al(t^{-1})$, and as $t$ increases from 0 to $\infty$,
so does $\al(t)$, monotonically.

\begin{lemma}
For $t > 0$, the sextic $\Phi_t(x,y):= \al(t)x^6 - x^4y^2 - x^2y^4
+ \be(t)y^6$ is positive with zeros at $(1,\pm t)$.
\end{lemma}
\begin{proof}
A  computation shows that 
\begin{equation}
\Phi_t(x,y) = \frac {(t^2x^2-y^2)^2((t^4+2t^2)x^2 + (2t^2+1)y^2)}{3t^4}.
\end{equation}
\end{proof}

Let $K= \{(x,y): x > 0,\  y \ge \ga(x))\}$ denote the region
lying above the curve $C = \{(\al(t),\be(t)): t > 0\}$, which
partially parametrizes the quartic curve $27x^2y^2
- 18x y - 4x - 4y - 1 =0$. For this reason,
\begin{equation}
\ga(x) = \frac{2+9x + 2(1+3x)^{3/2}}{27x^2}.
\end{equation}
\begin{lemma}
The binary sextic $\Psi(x,y) = r x^6 - x^4y^2 - x^2y^4 + s y^6$ is psd
if and only if $(r,s) \in K$. 
\end{lemma}
\begin{proof}
A necessary condition for the positivity of $\Psi$ is $r > 0$. Let $t_0
= \al^{-1}(r)>0$, so 
\begin{equation}
\Psi(x,y) = \Phi_{t_0}(x,y) + (s - \ga(t_0))y^6.
\end{equation}
If $(r,s) \in K$; that is, if $s\ge \ga(t_0)$, then Lemma 5.4 and
(5.8) show that 
$\Psi$ is positive. Conversely, $\Psi(1,t_0) =  (s - \ga(t_0))t_0^6$, 
so if $\Psi$ is positive, then $s\ge \ga(t_0)$.
\end{proof}

\begin{theorem}
The sextic $\Phi[c_1,0,c_3,1]$ is psd if and only if $(1+c_1,1+c_3)
\in K$.
\end{theorem}
\begin{proof}
One direction is clear by Lemmas 5.3(3) and 5.5. For the converse,
note that $(1+c_1,1+c_3) \in K$ if and only if $1+c_1 = \al(t_0)$
implies $1+c_3 \ge \be(t_0)$. In other words, we need to show that,
with $\la = 1 + c_3 - \be(t_0)$, 
\begin{equation}
\Phi[\al(t_0) -1,0,\be(t_0)+\la -1,1] = \Phi[\al(t_0) -1,0,\be(t_0)-1,1]
+ \la F_2^2
\end{equation}
is psd whenever $\la \ge 0$. To this end, for $t >0$, define
\begin{equation}
\begin{gathered}
R_t(x,y,z) := \Phi[\al(t) -1,0,\be(t)-1,1](x,y,z) = \\ 
\left(\frac{t^4+2t^2-3}{3} \right)^2 F_1^2(x,y,z)+
\left(\frac{1+2t^2-3t^4}{3t^4}\right)F_2^2(x,y,z) + R(x,y,z). 
\end{gathered}
\end{equation}
Note that $R_1 = R$, $R_{1/t}(x,y,z) = R_t(y,x,z)$ and that for $t \neq
1$, the coefficients of $F_1^2$ and $F_2^2$ have opposite sign.
 The following algebraic identity gives $Q_tR_t$ as
a sum of four squares for a psd  quadratic form $Q_t(x,y)$, which
implies that $R_t$ is psd, and completes the proof. 
\begin{equation}
\begin{gathered}
((2t^4+t^2)x^2+(t^2+2)y^2)3t^4R_t(x,y,z) 
\\= 3t^6(1 + 2t^2)x^2z^2(x^2 - z^2)^2 + 3t^4(2 + t^2)y^2z^2(y^2 - z^2)^2 \\+
   t^2(t^2 - 1)^2x^2y^2(t^2x^2 - y^2 + (1 - t^2)z^2)^2\\ + (2 + t^2)(1
   + 2t^2)(t^4x^4 - y^4 - t^4x^2z^2 + y^2z^2)^2. \end{gathered}
\end{equation}
\end{proof}
For $t=1$, (5.11) essentially appears in \cite[p.273]{Ro}.
In view of the foregoing, $\mz(R_t)$ contains, at least, ${\mathcal A} \cup
\{(1,\pm t,0)\}$. If $R_t(a,b,c)=0$, then each of the squares in (5.9)
vanishes.  In particular, $cF_1(a,b,c)=cF_2(a,b,c) = 0$, so either
$c=0$ or $(a,b,c) \in {\mathcal A} \cup \{(0,0,1)\}$. These cases have
already been 
discussed and we may conclude that $\mz(R_t)= {\mathcal A} \cup
\{(1,\pm t,0)\}$ and $|\mz(R_t)| = 10$. 

We now complete our discussion of the psd case.
\begin{theorem}
The sextic $\Phi[c_1,c_2,c_3,1]$ is psd if and only if $(c_1,c_3) \in
K$ and $|c_2| \le \sigma(c_1,c_3)$ for a function  $\sigma(c_1,c_3) \ge 0$
defined on $K$ (see (5.15)). If $c_2 = \pm \sigma(c_1,c_3)$, then
$\Phi[c_1,c_2,c_3,1]=R_t+\al(t^3 F_1 \pm F_2)^2$ (for suitable $t,
\al$ and choice of sign).
\end{theorem}

\begin{proof}
First, suppose  $\Phi[c_1,c_2,c_3,1]$ is psd. Then $(c_1,c_3) \in K$ by Lemma 
5.3(4) and Theorem 5.6. Setting $z=0$, we obtain the psd binary sextic 
\begin{equation}
\Gamma(x,y) = (1+c_1)x^6 - x^4y^2 + 2c_2 x^3y^3 - x^2y^4 + (1+c_3)y^6.
\end{equation}
Define $t_0$ so that $1+c_1 = \al(t_0)$. If $1+c_3 = \be(t_0)$, then
$\Gamma(1,\pm t_0) = \pm c_2t_0^3$ implies that $c_2 = 0$; otherwise,
$(1+c_1,1+c_3)$ lies 
strictly above $C$. Suppose now that $c_2 < 0$ without loss of
generality (taking $y \mapsto -y$ if necessary), so that for $u > 0$,
\begin{equation}
\Gamma(1,-u) > \Gamma(1,u) = (1+c_1) - u^2 - 2|c_2| u^3 - u^4
+(1+c_3)u^6 \ge 0.
\end{equation}
Let $\Psi(u) = (1+c_1)u^{-3} - u^{-1} - u + (1+c_3)u^3$, so that
\begin{equation}
0 \le u^{-3}\Gamma(1,u) = u^3(\Psi(u) - 2 |c_2|).
\end{equation}
Now define 
\begin{equation}
\sigma(c_1,c_3): = \min_{u > 0} \tfrac 12 \Psi(u) = \tfrac 12 \Psi(v);
\end{equation}
since $\Psi(u) \to \infty$ as $t \to 0$ or $t \to \infty$,
the minimum exists. It follows that $|c_2| \le
\sigma(c_1,c_3)$. (Although $\sigma(c_1,c_3)$ is computable explicitly, it is
quite complicated. For example, 
$2\sigma(1,0)$ is the unique real positive root of the sextic
$729x^6 - 22518 x^4 + 182774 x^2 - 111392$, approximately $.81392$.)

We must now show that every  $\Phi[c_1,\pm\sigma(c_1,c_3),c_3,1]$ is
psd. Since $\Psi'(v) = 0$, we have the system
\begin{equation}
\begin{gathered}
\sigma(c_1,c_3) = \frac12 \left( (1+c_3)v^3 - v - v^{-1} +
(1+c_1)v^{-3}\right); \\ 
3(1+c_3)v^2 - 1 + v^{-2} - 3(1+c_1)v^{-4} = 0.
\end{gathered}
\end{equation}
A calculation shows that (5.16) implies
\begin{equation}
\Phi[c_1,-\sigma(c_1,c_3),c_3,1] = R_v + \mu (v^3F_1 -F_2)^2,
\end{equation}
where $R_v$ is defined in (5.10) and
\begin{equation}
\mu = \frac{3(1+c_3)v^4 -(2v^2+1)}{3v^4}. 
\end{equation}
We are done if we can show that $\mu \ge 0$. 
By hypothesis, both sides of (5.17) vanish at $(1,v,0)$. But if we
evaluate (5.17) at $(1,-v,0)$, we have already seen that the left-hand side is
positive, and the right-hand side is  $0+4v^6\mu$, hence $\mu > 0$. 
\end{proof}

If $\Phi[c_1,c_2,c_3,1](a,b,c) = 0$, then Theorem 5.2 implies that
$(a,b,c) \in \mathcal A$ or
\begin{equation}
abc(a^2-c^2)(b^2-c^2)(a^2-ab+b^2-c^2)(a^2+ab+b^2-c^2)= 0.
\end{equation}
This includes the new zeros of $R_t$ on $c=0$ but also the extraneous
points $(a,b,c)$ for which $a^2+b^2-c^2 = \pm ab$, which never appear
non-trivially as zeros for any $R_t$.

To sum up, we have described sections of the two cones 
\begin{equation}
\begin{gathered}
P = \{(c_1,c_2,c_3,c_4) : c_1F_1^2 +2c_2F_1F_2 + c_3F_2^2 + c_4R  \in
P_{3,6} \} \subseteq \mathbb R^4, \\ 
\Sigma= \{(c_1,c_2,c_3,c_4) : c_1F_1^2 +2c_2F_1F_2 + c_3F_2^2 + c_4R
\in \Sigma_{3,6} \} \subseteq \mathbb R^4;
\end{gathered}
\end{equation}
at  $c_4=0$ and at $c_4=1$.  
 In the first case, the sections coincide and are literally a right
 regular cone. In the second case 
$\Sigma$ disappears, and if we think of $(c_1,c_3)$ as lying in a
plane and $c_2$ as the vertical dimension, then $P$ is a kind of
clam-shell, with a convex boundary curve $C$ lying in the plane and
rays emanating at 
varying angles from the points on the boundary.
 
\section{More ternary sextic examples}

\begin{example}
Let $A= \{\pi_i\} = \{(a_i,b_i)\}$ 
be given by $\pi_1 = (-1,0), \pi_2 = (-1,-1), \pi_3 = (0,1), \pi_4 =
(0,-1), \pi_5 = (1,0), \pi_6 = (2,2), \pi_7 = (2,-2), \pi_8 =
(1,-3)$. By looking at the $3 \times 3$ minors of the matrix with rows
$(1,a_i,b_i)$ and the $6 \times 6$ minors of the matrix with rows
$(1,a_i,b_i,a_i^2,a_ib_i,b_i^2)$, one can check that no three of the
$\pi_i$'s lie in a line, and no six on a quadratic. According to
Mathematica, $I_{1,3}(A)$ is spanned by
\begin{equation}
\begin{gathered}
f_1(x,y) = -42 + 49 x + 42x^2 - 49x^3 - 20 y - 38xy + 4x^2y + 42y^2 + 
    20y^3, \\
f_2(x,y) = -22 + 31 x + 22 x^2 - 31x^3 - 12y - 18xy + 22y^2 + 
    4xy^2 + 12y^3,
\end{gathered}
\end{equation}
and $\tilde A = \{\left(\frac{2516}{1297},\frac{4991}{2594}\right)\}$,
so $A$ is 
copacetic. In Hilbert's notation, $\phi(x,y) = x^2 -xy+y^2-1$ and 
\begin{equation}
\begin{gathered}
\psi(x,y) = -6136 + 2924x + 5784x^2 - 2924x^3 + 352x^4\\ - 2804y - 
    7000xy + 6299x^2y - 1049x^3y + 5818y^2\\ - 7803xy^2 + 
    1811x^2y^2 + 2804y^3 - 1402xy^3 + 318y^4.
\end{gathered}
\end{equation}
It follows that there exists $c>0$ so that $f_1^2+f_2^2 + c\phi\psi$
is psd and not sos. We do not offer an estimate for $c$.
\end{example}

In the examples in the rest of this section, the
symmetries are more clearly seen when the polynomials are homogenized. 

\begin{example}
We present one of several ways to generalize Robinson's
original set of eight points. For $t > 0$, let
\begin{equation}
A_t = \{(\pm 1,\pm 1), (\pm t, 0), (0, \pm t)\}.
\end{equation}
It is not hard to see that $A_t$ is copacetic (with ninth point
$(0,0)$) unless $t = \sqrt 2$, in 
which case $A_t$ lies on $x^2 + y^2 =2$. Since $A_t \mapsto
A_{2/t}$ under the invertible map
$(x,y) \mapsto ((x+y)/t,(x-y)/t)$, 
we may assume $0 < t < \sqrt 2$. After homogenizing to
$\mathcal A_t $, we note that a basis of $I_{1,3}({\mathcal A_t})$ is given by 
\begin{equation}
\{F_{1,t},F_{2,t}\} = 
\{x(x^2 + (t^2-1)y^2 - t^2z^2), y((t^2-1)x^2 + y^2 - t^2z^2)\}
\end{equation}
and that $\tilde{\mathcal A_t} = (0,0,1)$.
It is  not hard to see that
\begin{equation}
G_t(x,y,z) = 
(x^2 + (t^2-1)y^2 - t^2z^2)((t^2-1)x^2 + y^2 - t^2z^2)(-x^2-y^2+t^2z^2)
\end{equation}
is singular on $\mathcal{A}_t$ and is positive on $(0,0,1)$. 
(Robinson's example is recovered by setting $t=1$.) 

Consider now
\begin{equation}
\begin{gathered}
P_t := F_{1,t}^2 + F_{2,t}^2 + 1\cdot G_t^2 =  (2-t^2)(x^6
-x^4y^2-x^2y^4+ y^6) + \\ 
(2t^4 - 3t^2)(x^4+y^4)z^2 +(6t^2-4t^4+t^6)x^2y^2z^2 - t^6 (x^2z^4 +
y^2z^4-z^6).
\end{gathered}
\end{equation}
The proof that $P_t$ is psd follows from the identity
\begin{equation}
\begin{gathered}
(x^2+y^2)P_t = (2-t^2)(x^2-y^2)^2(x^2+y^2-t^2z^2)^2 + \\
t^2x^2z^2(x^2+(t^2-1)y^2-t^2z^2)^2
+t^2y^2z^2((t^2-1)x^2+y^2-t^2z^2)^2.
\end{gathered}
\end{equation}
For $t=1$, this formula is in \cite{Ro}. For $t = 0, \sqrt 2$, $P_t$
is sos. It is not hard to show that if $0<t<\sqrt 2$, then $\mz(P_t) =
{\mathcal A}_t \cup \{(1,\pm 1,0) \}$ has 10 points and $P_t$ is not sos.
\end{example}

\begin{example}
Let
\begin{equation}
{\mathcal A}  = \{(1,0,0),(0,1,0),(0,0,1),(1,1,0),(1,0,1),(1,1,0),(1,1,1) \}.
\end{equation}
It is again simple to show that $I_{1,3}({\mathcal A})$ is spanned by
\begin{equation}
F_1(x,y,z) = xy(x-y),\quad F_2(x,y,z) = yz(y-z),\quad  F_3(x,y,z) = zx(z-x),
\end{equation}
and that
\begin{equation}
G(x,y,z) = xyz(x-y)(y-z)(z-x)
\end{equation}
is in  $I_{2,6}({\mathcal A}) \smallsetminus I_{1,3}^2({\mathcal
  A})$. Accordingly, by Theorem 4.3, there exists $c>0$ so that
\begin{equation}
U_c(x,y,z) = x^2y^2(x-y)^2 + y^2z^2(y-z)^2 + z^2x^2(z-x)^2 +
cxyz(x-y)(y-z)(z-x) 
\end{equation}
is psd and not sos.
Since  $U_c(x,y,z) \ge 0$ whenever $xyz=0$, we define
\begin{equation}
\begin{gathered}
Q_c(x,y,z) := \frac{U_c(x,y,z)}{x^2y^2z^2} \\ = \frac {(x-y)^2}{z^2} +
\frac{(y-z)^2}{x^2} + \frac{(z-x)^2}{y^2} + c \left(\frac  {x-y}{z}\right)
 \left(\frac  {y-z}{x}\right) \left(\frac  {z-y}{x}\right).
\end{gathered}
\end{equation}
It is now sensible to make a substitution: let
\begin{equation}
u  := \frac  {x-y}{z};\quad v := \frac  {y-z}{x}; \quad w  := \frac  {z-x}{y}.
\end{equation}
Then $Q_c = u^2 + v^2 + w^2 + cuvw$; somewhat surprisingly,
$\{u,v,w\}$ is not algebraically independent: in fact,
\begin{equation}
u+v+w+uvw = 0.
\end{equation}
An application of Lagrange multipliers to minimize $Q_c$, subject to
(6.14), shows that two of $\{u,v,w\}$ are equal; by symmetry, we
may take $u=v$, so that $w = -\frac{2u}{u^2+1}$, and
\begin{equation}
Q_c\left(u, u, -\tfrac{2u}{u^2+1}\right)  =
\frac{2u^2(u^4+2u^2+3 -cu(1+u^2))}{(1+u^2)^2}.
\end{equation}
Let $\sigma = \sqrt{\sqrt 2 + 1}$. A little calculus
shows that the numerator is psd provided $|c| \le c_0:= 
4/\sigma$, with $Q_{c_0} = 0$ when $u = \pm \sigma$. Solving back for
$(x,y,z)$ yields, up to multiple, that $(1+\sigma,
1+\sigma^2,1-\sigma)$ and its cyclic images are in $\mz(U_{c_0})$, together
with (6.8). Here, $| \mz(U_{c_0}) | = 10$.
\end{example}
 
\begin{example}
The Motzkin form $M$ cannot be derived directly from Theorems 4.1 or
4.3  because $|\mz(M)| = 6$;
however, $M$ has zeros at $(1,0,0)$ and $(0,1,0)$ which vanish
to the sixth order in the $z$-direction. It is possible to construct
psd ternary sextics $M_t$ with $|\mz(M_t)| = 10$ for $t>0$ and such
that $M_t \to M$ as $t \to 0$. We do this with an Ansatz by supposing
that there is a non-zero even ternary sextic which is
 symmetric in $(x,y)$ and lies
in $I_{2,6}({\mathcal A_t})$ for
\begin{equation}
{\mathcal A_t} = \{(1,0,0), (0,1,0), (1, 0, \pm t), (0, 1, \pm t), (1,
\pm 1, \pm 1)\}. 
\end{equation} 
Although these impose 30 equations on the 28 coefficients of a ternary
sextic, there is some redundancy, and it can be verified that
\begin{equation}
\begin{gathered}
M_t(x,y,z) = (1-2t^2)(x^4y^2+x^2y^4) + t^4(x^4z^2+y^4z^2)\\- (3 -
8t^2+2t^4)x^2y^2z^2  -2t^2(x^2+y^2)z^4 + z^6
\end{gathered}
\end{equation}
satisfies this criterion. It is
not clear  that $M_t$ is psd; in fact, it is not psd when
$t^2 > 1/2$. We 
note that $M_0 = M$ and $M_t$ is a square when $t^2 = 1/2$. The
proof that $M_t$ is psd for $t^2 < 1/2$ is given by an sos representation of
$Q_tM_t$:
\begin{equation}
\begin{gathered}
(x^2+y^2)M_t(x,y,z) = (1-2t^2)x^2y^2(x^2+y^2-2z^2)^2 + \\
  y^2z^2(t^2(x^2-y^2)-(x^2-z^2))^2 + x^2z^2(t^2(y^2-x^2)-(y^2-z^2))^2.
\end{gathered}
\end{equation}
This equation also shows that, at least when $t^2 < 1/2$,  $\mz(M_t)=
\mathcal A_t$.
We may also derive $M_t$ using Theorem 4.1, by first
choosing any eight points in $\mathcal A_t$.
\end{example}

\begin{example}
Similarly, one can approach $S(x,y,z)$ by Ansatz and look for a
cyclically symmetric even sextic $S_t$ which is singular at 
\begin{equation}
{\mathcal A_t} = \{(\pm t,1,0), (0, \pm t,1), (1,0,\pm t), (1, \pm 1,
\pm 1)\}. 
\end{equation}
Again, although there is no reason to expect a non-zero solution,
there is one:
\begin{equation}
\begin{gathered}
S_t(x,y,z) = t^4(x^6+y^6+z^6) + (1-2t^6)(x^4y^2+
y^4z^2+z^4x^2)\\ +
(t^8 - 2t^2)(x^2y^4+y^2z^4+z^2x^4)-3(1-2t^2+t^4-2t^6+t^8)x^2y^2z^2.
\end{gathered}
\end{equation}
We find that
$t^8S_{1/t}(x,y,z) = S_t(x,z,y)$,  $S_0(x,y,z) = S(x,y,z)$ and
$S_1(x,y,z) = R(x,y,z)$. 
 The proof that $S_t$ is psd follows from yet
another algebraic identity: 
\begin{equation}
\begin{gathered}
(x^2+y^2)S_t(x,y,z) =
(t^2x^4+x^2y^2-t^4x^2y^2-t^2y^4-x^2z^2+t^4y^2z^2)^2\\ +
 y^2z^2(y^2-x^2+t^2(x^2-z^2))^2 + t^4x^2z^2(y^2-z^2 +
 t^2(x^2-y^2))^2 \\
+(t^2-1)^2x^2y^2((z^2-x^2)+t^2(y^2-z^2))^2.
\end{gathered}
\end{equation}
When $t=1$, (5.11) and (6.21) coincide. This example was
announced, without proof, in \cite[p.261]{Re2}.
\end{example}

Robinson \cite[p.273]{Ro} observed that $(ax^2+by^2+cz^2)R(x,y,z)$ is sos, ``at
least if $0 \le a \le b+c,\ 0 \le b \le a+c,\ 0 \le c \le a+b$."
We revisit this situation and simultaneously
illustrate the method used to discover (5.11), (6.7), (6.18) and (6.21). 

\begin{theorem}
If $r,s,t \ge 0$, then
$(r^2x^2 + s^2y^2 + t^2z^2)R(x,y,z)$ is sos if
and only if $r \le s+t$,  $s \le r+t$ and $t \le r+s$.
\end{theorem}
\begin{proof}
It was shown in \cite[p.569]{CLR2} (by a polarization argument)
that an even sos polynomial  $F$ has an sos
representation $F = \sum H_j^2$ in which each $H_j^2$ is even. Suppose
\begin{equation}
(r^2x^2 + s^2y^2 + t^2z^2)R(x,y,z) = \sum_{j=1}^r H_j^2(x,y,z)
\end{equation}
is such an ``even'' representation.
Then  $\mz(R) \subseteq \mz(H_j)$ for the quartic $H_j$'s
(c.f. (5.4)).  It follows that
\begin{equation}
\begin{gathered}
H_j(x,y,z) = c_{1j} xy(x^2-y^2) + c_{2j} xz(x^2-z^2) + c_{3j}
yz(y^2-z^2) \\ + (c_{4j}(x^2-z^2)(x^2-y^2+z^2)+c_{5j}(y^2-z^2)(-x^2+y^2+z^2)). 
\end{gathered}
\end{equation}
Each $H_j^2$ is even, so the only 
cross-terms which can appear in any $H_j^2$ are $c_{4j}c_{5j}$ and 
\begin{equation}
\begin{gathered}
(r^2x^2 + s^2y^2 + t^2z^2)R(x,y,z) = \la_1x^2y^2(x^2-y^2)^2 +
  \la_2x^2z^2(x^2-z^2)^2 \\ + \la_3 y^2z^2(y^2-z^2)^2    + \la_4
  (x^2-z^2)^2(x^2-y^2+z^2)^2 +\\ 2\la_5
  (x^2-z^2)(x^2-y^2+z^2)(y^2-z^2)(-x^2+y^2+z^2) \\+ \la_6
  (y^2-z^2)^2(-x^2+y^2+z^2)^2,
\end{gathered}
\end{equation}
for $\la_j$'s, defined by
\begin{equation}
\begin{gathered}
\la_1 = \sum_j c_{1j}^2,\quad \la_2 = \sum_j c_{2j}^2, \quad \la_3 = \sum_j
c_{3j}^2,\\ \la_4 = \sum_j c_{4j}^2,\quad  \la_5 = \sum_j
c_{4j}c_{5j}, \quad  \la_6 = \sum_j c_{5j}^2.
\end{gathered}
\end{equation}
We solve for the $\la_j$  in (6.24):
\begin{equation}
\la_1 = t^2,\quad \la_2 = s^2,\quad \la_3 = r^2,\quad \la_4 =
r^2,\quad \la_6 = s^2,\quad \la_5 = (t^2-r^2-s^2)/2.  
\end{equation}
There exist $c_{ij}$ to satisfy (6.25) and (6.26) if and only if
\begin{equation}
0 \le \la_4\la_6- \la_5^2 = \frac14 (r+s-t)(r+t-s)(s+t-r)(r+s+t) 
\end{equation}
If, say, $r \ge s \ge t\ge 0$, then $r+s\ge t$ and  $r+t \ge s$
automatically, and so (6.27) holds if and only if $s+t \ge r$. By
symmetry, we see that (6.27) is true if and only if all three
inequalities hold.
\end{proof}

\section{Extremal psd ternary forms}

In 1980, Choi, Lam and the author \cite{CLR1} studied
$|\mz (F)|$ for $F \in P_{3,m}$. Let 
\begin{equation}
\al(m):= \max \left( \frac {m^2}4, \frac {(m-1)(m-2)}2 \right).
\end{equation}
By Theorem 3.5 in \cite{CLR1}, if $F \in P_{3,m}$, then $|\mz(F)| >
\al(m)$ implies $|\mz(F)| = 
\infty$, and this occurs if and only if $F$ is divisible by the square
of an indefinite form. Let
\begin{equation}
B_{3,m} = \{ \sup |\mz(F)| \ : \ F \in P_{3,m},\  |\mz(F)| < \infty \}. 
\end{equation}
Then by Theorem 4.3 in \cite{CLR1},
\begin{equation}
\begin{gathered}
\frac{m^2}4 \le B_{3,m} \le \frac {(m-1)(m-2)}2; \\B_{3,6k} \ge 10k^2,\quad
B_{3,6k+2} \ge 10k^2+1,\quad B_{3,6k+4} \ge 10k^2+4.
\end{gathered}
\end{equation}
In particular, $B_{3,6} = 10$. Further, if  $F \in P_{3,6}$, and $|\mz(F)| >
10$, then $|\mz(F)| = \infty$ and
$F\in \Sigma_{3,6}$ is a sum of three squares (Theorem 3.7).
If $G$ is a ternary sextic and   
$|\mz(G)| = 10$, then one of $\pm G$ is psd and not sos (Corollary
4.8). We wrote (p.12): ``it would be of interest to determine, if
possible, {\it all} forms $p \in P_{3,6}$ with exactly 10 zeros. From
a combinatorial point of view, it would already be of interest to
determine (or classify) all configurations of 10-point sets $S \subset
\mathbb P^2$ for which there exist $p \in P_{3,6}$ such that $S =
\mz(p)$ $\dots$ The only known psd ternary sextic
with 10 zeros is $R$.''
Sections five and six of this paper are inspired by this remark.
\begin{lemma}
If $F \in P_{3,6}$ is reducible, then $F \in \Sigma_{3,6}$.
\end{lemma}
\begin{proof}
If $F$ has an indefinite factor $H$, then $F = H^2G$, where
$G \in P_{3,2d} = \Sigma_{3,2d}$ for $2d \le 4$. If $F = F_1F_2$ for
definite $F_i$, then $\deg F_i \le 4$ again implies $F \in \Sigma_{3,6}$.
\end{proof}

A form $F$ in the closed convex cone 
$P_{n,m}$ is {\it extremal} if $F = G_1 + G_2$ for $G_j \in
P_{n,m}$ implies that $G_j = \la_j F$ for $0 \le \la_j \in
\rr$. Equivalently, $F$ is extremal if $F \ge G \ge 0$ implies $G =
\la F$. The set of extremal forms in $P_{n,m}$ is  denoted by
$E(P_{n,m})$.

\begin{theorem}
Suppose $F \in P_{3,6}$ and $|\mz(F)| = 10$. Then $F \in E(P_{3,6})$. 
\end{theorem}
\begin{proof}
Since $F \in \Delta_{3,6}$ by \cite{CLR1}, Lemma 7.1 implies that $F$ is
irreducible. Suppose $F \ge G \ge 0$. Then $F$ and $G$ are both
singular at the ten zeros of $F$, and since $10\cdot 2^2 > 6\cdot 6$,
Bezout implies that $F$ and $G$ have a common factor. Thus $G = \la F$
and $F$ is extremal. 
\end{proof}

Theorems 5.1 and 5.7 imply that if $F \in E(P_{3,6})$ has Robinson's 8
zeros, then either $F = P_t \in \Delta_{3,6}$ for some $t > 0$ has ten
zeros, or $F = (\al F_1 + \be F_2)^2 \in E(\Sigma_{3,6})$.

We can use the Perturbation Lemma to put a strong restriction on those 
extremal forms which only have round zeros. 

\begin{theorem}
If $P \in E(P_{3,2d})\cap \Delta_{3,2d}$ and all zeros of $P$ are
round, then $|\mz(P)| \ge \frac{(d+1)(d+2)}2$. 
\end{theorem}
\begin{proof}
Suppose $P$ is psd, all its zeros are round, and $|\mz(P)|<
\frac{(d+1)(d+2)}2$. Then there exists a non-zero $H \in
I_{1,d}(\mz(P))$
and the Perturbation Lemma applies to $(P,\pm H^2)$.
 It follows that $P \pm cH^2$ is psd for some $c>0$ and
$P$ is not extremal because
\begin{equation}
P= \tfrac 12(P - cH^2) + \tfrac 12(P+ cH^2);
\end{equation}
$P \neq \la H^2$ since $P$ is not sos.
\end{proof}
\begin{corollary}
If $p \in E(P_{3,6})\cap \Delta_{3,6}$ and all zeros of $P$ are
round, then $|\mz(p)|=10$.
\end{corollary}

\begin{lemma}
If $P \in P_{3,6}$, and $\mz(P)$ contains four points in a line or
seven points on a quadratic, then $P \in \Sigma_{3,6}$.
\end{lemma}
\begin{proof}
If $\mz(P)$ contains four points $\pi_i$ on the line $L$, then
since $P$ is singular at its zeros, Bezout implies that
$L$ divides $P$ and $P \in \Sigma_{3,6}$ by Lemma 7.1. Similarly, if
$\mz(P)$ contains seven points $\pi_i$ on the quadratic $Q$, then 
Bezout again implies that $P$ is reducible.
\end{proof}

\begin{theorem}
If $P \in E(P_{3,6})\cap \Delta_{3,6}$ and all zeros of $P$ are
round, then $P$ can be derived by Hilbert's Method using Theorem 4.3.
\end{theorem}
\begin{proof}
Let $A$ denote any subset of seven of the ten zeros of $P$. By Lemma
7.5, $A$ meets the hypothesis of Theorem 4.3.
\end{proof}

Given positive $f \in \mathbb R_{n,2d}$ and $\pi \in \mathbb R^n$, let
$E(f,\pi)$ denote the set of $g \in \mathbb R_{n,d}$  
such that there exists a neighborhood ${\mathcal N}_g$ of $\pi$ and $c
> 0$ so that $f - cg^2$ is non-negative on  
${\mathcal N}_g$.

\begin{lemma}
$E(f,\pi)$ is a subspace of $\mathbb R_{n,d}$. 
\end{lemma}
\begin{proof}
Clearly, $g \in E(f,\pi)$ implies $\la g \in
E(f,\pi)$ for $\la \in \mathbb R$. Suppose $g_1, g_2 \in E(f,\pi)$;
specifically, $f - c_1g_1^2 \ge 0$ on ${\mathcal N}_1$ and  $f -
c_2g_2^2 \ge 0$ on ${\mathcal N}_2$, and let  
$\mathcal N = {\mathcal N}_{1} \cap {\mathcal N}_{2}$ and $c =
\min(c_1,c_2)$. The identity
\begin{equation}
f - \tfrac c4(g_1+g_2)^2 = \tfrac 12( f - cg_1^2) + \tfrac 12(f -
cg_2^2) + \tfrac c4(g_1-g_2)^2 
\end{equation}
shows that $g_1+g_2 \in E(f,\pi)$.
\end{proof}
 If $f(\pi) > 0$, then $E(f,\pi) = \mathbb R_{n,d}$. Let 
\begin{equation}
\de(f,\pi) := \binom{n+d}d - \dim E(f,\pi)
\end{equation}
measure the singularity of the zero of $f$ at $\pi$; the argument of
the Perturbation Lemma shows that
$\de(f,\pi) = 1$ if and only if $f$ has a round zero at $\pi$.  These
definitions also apply in the obvious way to the homogeneous case.

\begin{theorem}
If $P \in E(P_{3,2d})\cap\Delta_{3,2d}$, then
\begin{equation}
\delta(P): = \sum_{\pi \in \mz(P)} \de(P,\pi) \ge \frac{(d+1)(d+2)}2.
\end{equation}
\end{theorem}
\begin{proof}
 If $f(\pi) > 0$, then $E(f,\pi) = \mathbb
R_{n,d}$. 
Let 
\begin{equation}
{\mathcal E}:= \bigcap_{\pi \in \mz(P)} E(f,\pi).
\end{equation}
Since
\begin{equation}
\dim {\mathcal E} \ge \frac{(d+1)(d+2)}2 - \de(P),
\end{equation}
if (7.7) fails, then there exists $0 \neq H \in {\mathcal E}$. The
argument of 
the Perturbation Lemma applies to $(P,\pm H^2)$, so that (7.4)
holds for some $c > 0$, and $P$ is not extremal. 
\end{proof}

It can be checked that $M$ has round zeros at $(1,\pm 1,
\pm 1)$. Let $\pi =  (1,0,0)$. If $M - c F^2$ is non-negative near
$(1,0,0)$ for a ternary cubic $F$, then by the method of cages (see
\cite[\S 3]{CLR3}), $x^3, x^2z, xz^2$ cannot appear in $F$, whereas every
other monomial is in $E(M,\pi)$, and so $\delta(M,\pi) =
3$. By symmetry,  $\delta(M,(0,1,0)) = 3$,
so that $\de(M) =  4\cdot 1 + 2 \cdot 3 = 10.$
A similar calculation for $S$ shows that it has round zeros at $(1,\pm 1,
\pm 1)$ and that  $\delta(S,e_i) = 2$ at the unit vectors $e_i$
so $\de(S) = 4 \cdot 1 + 3 \cdot 2 = 10$ as well. 
Examples 6.4 and 6.5 were constructed under a heuristic in which
``coalescing'' zeros explain higher-order singularities. These lead to
a perhaps overly-optimistic conjecture:

\begin{conjecture}
If $P \in E(P_{3,6}) \cap \Delta_{3,6}$, then $\de(P) = 10$,
and either $P$ has ten round zeros, or is the limit of psd extremal
ternary sextics with ten round zeros.
\end{conjecture}

These results are likely more complicated in higher
degree.  The ternary octic
\begin{equation}
T(x,y,z) = x^4y^4 + x^2z^6+y^2z^6 - 3x^2y^2z^4 = x^4y^4z^6M(1/x,1/y,1/z)
\end{equation}
is in $E(P_{3,8}) \cap \Delta_{3,8}$; see \cite[p.372]{Re0}. It has
five round zeros at $(0,0,1)$ and $(1,\pm1, \pm 1)$, and more
singular zeros at $(1,0,0)$ and $(0,1,0)$ at which $\delta = 5$,
so that $\delta(T) = 15$. On the other hand, for 
\begin{equation}
U(x,y,z) = x^2(x-z)^2(x-2z)^2(x-3z)^2 + y^2(y-z)^2(y-2z)^2(y-3z)^2 \in
\Sigma_{3,8},
\end{equation}
$\mz(U) = \{(i,j,1): 0 \le i, j \le 3\}$, so $\delta(U) = 16$. 
Thus, there is no threshold value for $\delta$ separating
$\Sigma_{3,8}$ and $\Delta_{3,8}$, as there is for sextics.
\section{Ternary forms in higher degree}
For $d \ge 3$, let
\begin{equation}
T_d = \{(i,j)\ : \ 0 \le i, j,\ i+j \le d\} \subset \mathbb Z^2
\end{equation}
denote  a right triangle of $\frac{(d+1)(d+2)}2$ lattice points. 
Define the falling product by
\begin{equation}
(t)_m=\prod_{j=0}^{m-1}(t-j).
\end{equation}

The following
construction is due to Biermann \cite{Bie}, see \cite[pp.31-32]{Re1}.
For $(r,s) \in T_d$, let
\begin{equation}
\phi_{r,s,d}(x,y) := \frac{ (x)_r(y)_s(d-x-y)_{d-r-s}}{r!s!(d-r-s)!}.
\end{equation}
\begin{lemma}
If $(i,j) \in T_d$, then $\phi_{r,s,d}(i,j) = 0$ if $(i,j) \neq (r,s)$
and $\phi_{r,s,d}(r,s) = 1$.
\end{lemma}
\begin{proof}
Observe that $(n)_m=0$ if $n \in \{0,\dots,m-1\}$ and $(m)_m =
m!$. If $(i,j) \in T_d$, then $0 \le i$, $0 \le j$ and $0 \le d - i -
j$. Thus $\phi_{r,s,d}(i,j) = 0$ unless $i \ge r$, $j \ge s$ and $d-i-j
\ge d - r - s$, or $i+j \le r+s$; that is, unless $(i,j) = (r,s)$.
The second assertion is immediate.
\end{proof}

\begin{theorem}
Suppose $B \subseteq T_d$ and $A = T_d \smallsetminus B$. Then a basis for
$I_{1,d}(A)$ is given by $\{\phi_{r,s,d} : (r,s) \in B\}$.
\end{theorem}

\begin{proof}
The set $\{\phi_{r,s,d}: (r,s) \in T_d\}$ consists of the correct
number of linearly independent polynomials and so is a basis for
$\rr_{2,d}$.
If $p \in  \rr_{2,d}$, then upon evaluation at $(r,s)
\in T_d$, we immediately obtain
\begin{equation}
p(x,y) = \sum_{(r,s) \in T_d} p(r,s)\phi_{r,s,d}(x,y).
\end{equation}
If $p \in I_{1,d}(A)$, then $\phi_{r,s,d}$ has non-zero
coefficient in (8.4) only if $(r,s) \in B$.
\end{proof}

We use this construction in the following example, which was inspired by
looking at the regular pattern of pine trees below the Sulphur
Mountain tram, during a break in the October 2006 BIRS program on
``Positive Polynomials and Optimization''.
\begin{example}[The Banff Gondola Polynomials]
Suppose $d \ge 3$ and let
\begin{equation}
A_d = T_d \smallsetminus \{(d,0),(0,d)\} = \{(i,j): 0 \le i,j \le d-1,
i+j \le d\}.
\end{equation}
By Theorem 8.2, $I_{1,d}(A_d)$ is spanned by $f_1(x,y) =
\phi_{d,0,d}(x,y) = (x)_d$ and 
$f_2(x,y) = \phi_{0,d,d}(x,y) = (y)_d$, and it is easy to see that
$\mz(f_1) \cap \mz(f_2) 
= \{0,\dots,d-1\}^2$, so that
\begin{equation}
\tilde A_d = \{(i,j): 0 \le i,j \le d-1, i+j \ge d+1\}.
\end{equation}
Note that $(i,j) \in \tilde A_d$ implies that $i,j \ge 2$.
 Let 
\begin{equation}
\begin{gathered}
g_d(x,y) = (x)_2(y)_2(x+y-2)_{d-1}(x+y-4)_{d-3} \\=
x(x-1)y(y-1)(x+y-2)(x+y-3) \prod_{k=0}^{d-4} (x+y-4-k)^2,
\end{gathered}
\end{equation}
We claim that $g_d$ is singular at $\pi \in A_d$ and positive at $\pi \in
\tilde A_d$. First, it is easy to check that each point
in $A_3$ lies on at least two of the lines, and $g_3(2,2) =8$. Now
suppose $d \ge 4$ and $(r,s) \in A_d$. If $4 \le r+s \le d$, then
$(r,s)$ lies on a squared factor; if $2 \le r+s \le 3$, then $(r,s)$
lies on $x+y-2=0$ or $x+y-3=0$, but also, at least one of $\{r,s\}$ is 0
or 1. Finally, if $0 \le r+s \le 1$, then $\{r,s\}\subseteq \{0,1\}$.
If $(r,s) \in \tilde A_d$ for any $d$, then $r,s \ge 2$ and $r+s \ge d+1$,
so each factor in $g_d$ is positive at $(r,s)$. It follows from
Theorem 3.4 that there exists $c_d > 0$ so that
\begin{equation} 
(x)_d^2 + (y)_d^2 + c_d(x)_2(y)_2(x+y-2)_{d-1}(x+y-4)_{d-3}
\end{equation}
is positive and not a sum of squares. Note that this polynomial has at
least $|A_d|$ zeros, so $B_{3,2d} \ge \frac{d^2+3d-2}2$. This improves
the lower bound in (7.3) for $2d = 8, 10$.
It can be shown that $c(3) = 4/3$ (exactly) and
that $c(d) \le 12d^{-2}$, so $c(d) \to 0$. 
\end{example}
We conclude with some speculations about Hilbert's Method in
degree $d\ge 4$. Suppose $A$ is a set of $\binom {d+2}2-2$
points in general position, so that $I_{1,d}(A)$ has basis
$\{f_1,f_2\}$. By Bezout, we can only say that $|\tilde A| \le d^2 - |A| =
\binom{d-1}2$ as the common zeros do not have to be real or distinct. We have
$\dim I_{1,d}^2(A) = 3$ and, from (2.5), 
\begin{equation}
\dim I_{2,2d}(A) \ge \binom{2d+2}2 - 3\left(\binom {d+2}2-2\right) =
\binom{d-1}2 + 3.
\end{equation}
There exist $\binom{d-1}2$ linearly independent polynomials
in $I_{2,2d}(A) \smallsetminus I_{1,d}^2(A)$, and it is plausible that 
 one  is positive on $\tilde A$. If so, then
 Hilbert's Method could be applied. 

If $r \ge 3$, and  $A$ is a set of $\binom {d+2}2-r$
points in general position, so that $\dim I_{1,d}(A) = r$, 
then it is plausible to expect $\tilde {\mathcal A} =
\emptyset$. We have
\begin{equation}
\begin{gathered}
\dim I_{2,2d}(A) \ge \binom{2d+2}2 - 3\left(\binom {d+2}2-r\right) =
\binom{d-1}2 + 3r-3\\ = \frac {r(r+1)}2 + \frac{(d+1-r)(d+r-4)}2
\ge \dim I_{1,d}^2 + \frac{(d+1-r)(d+r-4)}2,
\end{gathered}
\end{equation}
so if $r \le d$,  $I_{2,2d}(A) \smallsetminus
 I_{1,d}^2(A)$ would be non-empty, and again Hilbert's Method
could be applied. We hope to return to these questions elsewhere.


\end{document}